\documentclass[a4paper]{article}
\usepackage[utf8]{inputenc}
\usepackage[T1]{fontenc}
\usepackage{amsmath,amsthm,amssymb,amsfonts}
\usepackage{thmtools}
\usepackage{xspace}
\usepackage{hyperref}
\usepackage{verbatim}
\usepackage{fullpage}
\usepackage{tikz}
\usetikzlibrary{arrows.meta}
\usepackage[procnumbered,linesnumbered,ruled,vlined]{algorithm2e}
\usepackage{tabularray}
\usepackage[textsize=footnotesize,color=green!40]{todonotes}
\usepackage{float}
\usepackage{stmaryrd}
\usepackage{subcaption}
\usepackage{enumitem}
\usepackage{cleveref}

\newcommand{\annPath}{concerning path\xspace}
\newcommand{\annPaths}{concerning paths\xspace}
\newcommand{\badPath}{unfixable path\xspace}
\newcommand{\badPaths}{unfixable paths\xspace}
\newcommand{\notSoBadPath}{fixable path\xspace}
\newcommand{\notSoBadPaths}{fixable paths\xspace}

\newcommand{\MD}{\textsc{Metric Dimension}\xspace}

\newcommand{\nicearrow}{-{Latex[scale=1.1]}}

\DeclareMathOperator{\dist}{dist}

\newcommand{\mw}{{\rm mw}}
\newcommand{\vp}{\mathbf{p}}

\usetikzlibrary{positioning}
\usetikzlibrary{fit}
\usetikzlibrary{decorations.pathmorphing}
\usetikzlibrary{calc}
\usetikzlibrary{shapes}

\tikzstyle{noeud}=[circle,inner sep=2, minimum size =3 pt, line width = 1pt, draw=black, fill=white]
\tikzstyle{bignoeud}=[circle, draw, inner sep=2, minimum size=10pt, line width=1pt]
\newcommand{\convexpath}[2]{
	[   
	create hullnodes/.code={
		\global\edef\namelist{#1}
		\foreach [count=\counter] \nodename in \namelist {
			\global\edef\numberofnodes{\counter}
			\node at (\nodename) [draw=none,name=hullnode\counter] {};
		}
		\node at (hullnode\numberofnodes) [name=hullnode0,draw=none] {};
		\pgfmathtruncatemacro\lastnumber{\numberofnodes+1}
		\node at (hullnode1) [name=hullnode\lastnumber,draw=none] {};
	},
	create hullnodes
	]
	($(hullnode1)!#2!-90:(hullnode0)$)
	\foreach [
	evaluate=\currentnode as \previousnode using \currentnode-1,
	evaluate=\currentnode as \nextnode using \currentnode+1
	] \currentnode in {1,...,\numberofnodes} {
		-- ($(hullnode\currentnode)!#2!-90:(hullnode\previousnode)$)
		let \p1 = ($(hullnode\currentnode)!#2!-90:(hullnode\previousnode) - (hullnode\currentnode)$),
		\n1 = {atan2(\y1,\x1)},
		\p2 = ($(hullnode\currentnode)!#2!90:(hullnode\nextnode) - (hullnode\currentnode)$),
		\n2 = {atan2(\y2,\x2)},
		\n{delta} = {-Mod(\n1-\n2,360)}
		in 
		{arc [start angle=\n1, delta angle=\n{delta}, radius=#2]}
	}
	-- cycle
}

\definecolor{nicered}{RGB}{162, 27, 43}
\definecolor{nicegreen}{rgb}{0.0, 0.42, 0.24}

\newtheorem{theorem}{Theorem}
\newtheorem{definition}[theorem]{Definition}
\newtheorem{proposition}[theorem]{Proposition}

\Crefname{definition}{Definition}{Definitions}

\title{Algorithms and hardness for Metric Dimension on digraphs\footnote{Work financed by the French government IDEX-ISITE initiative CAP 20-25 (ANR-16-IDEX-0001), the International Research Center ``Innovation Transportation and Production Systems'' of the I-SITE CAP 20-25, and by the ANR project GRALMECO (ANR-21-CE48-0004).}}
\author{Antoine Dailly\footnote{Université Clermont Auvergne, CNRS, Clermont Auvergne INP, Mines Saint-\'Etienne, LIMOS, 63000 Clermont-Ferrand, France.}
	\and Florent Foucaud\footnotemark[2]
	\and Anni Hakanen\footnotemark[2] \footnote{Department of Mathematics and Statistics, University of Turku, FI-20014, Finland} \footnote{Turku Collegium for Science, Medicine and Technology TCSMT, Finland} \footnote{Research supported by the Jenny and Antti Wihuri Foundation and partially by Research Council of Finland grant number 338797.}
}

\begin{document}
\maketitle

\begin{abstract}
	In the \MD problem, one asks for a minimum-size set $R$ of vertices such that for any pair of vertices of the graph, there is a vertex from $R$ whose two distances to the vertices of the pair are distinct. This problem has mainly been studied on undirected graphs and has gained a lot of attention in the recent years. We focus on directed graphs, and show how to solve the problem in linear time on digraphs whose underlying undirected graph (ignoring multiple edges) is a tree. This (non-trivially) extends a previous algorithm for oriented trees. We then extend the method to orientations of unicyclic graphs. We also give a fixed-parameter-tractable algorithm for digraphs when parameterized by the directed modular-width, extending a known result for undirected graphs. Finally, we show that \MD is NP-hard even on planar triangle-free acyclic digraphs of maximum degree~6.
\end{abstract}

\section{Introduction}

The metric dimension of a (di)graph $G$ is the smallest size of a set of vertices that distinguishes all vertices of $G$ by their vectors of distances from the vertices of the set. This concept was introduced in the 1970s by Harary and Melter~\cite{HM76} and by Slater~\cite{S75} independently. Due to its interesting nature and numerous applications (such as robot navigation~\cite{KRR96}, detection in sensor networks~\cite{S75} or image processing~\cite{MT84}, to name a few), it has enjoyed a lot of attention. It also has been studied in the more general setting of metric spaces~\cite{B53}, and is generally part of the rich area of identification problems of graphs and other discrete structures~\cite{lobsteinBIB}.

More formally, let us denote by $\dist(x,y)$ the distance from $x$ to $y$ in a digraph. Here, the distance $\dist(x,y)$ is taken as the length of a shortest directed path from $x$ to $y$; if no such path exists, $\dist(x,y)$ is infinite, and we say that $y$ is not \emph{reachable} from $x$. We say that a set $S$ is a \emph{resolving set} of a digraph $G$ if for any pair of distinct vertices $v,w$ from $G$, there is a vertex $x$ in $S$ with $\dist(x,v)\neq \dist(x,w)$. Furthermore, we require that every vertex of $G$ is reachable from at least one vertex of $S$. The \emph{metric dimension} of $G$ is the smallest size of a resolving set of $G$, and a minimum-size resolving set of $G$ is called a \emph{metric basis} of $G$.\footnote{The definition that we use has been called \emph{strong metric dimension} in~\cite{ABCHONS20}, as opposed to \emph{weak metric dimension}, where one single vertex may be unreachable from any resolving set vertex. The former definition seems more natural to us. However, the term \emph{strong metric dimension} is already used for a different concept, see~\cite{strongMD}. Thus, to prevent confusion, we avoid the prefix \emph{strong} in this paper.}

We denote by \MD the computational version of the problem: given a (di)graph $G$, determine its metric dimension. 

For undirected graphs, \MD has been extensively studied, and its non-local nature makes it highly non-trivial from an algorithmic point of view. On the hardness side, \MD was shown to be NP-hard for planar graphs of bounded degree~\cite{MDplanar}, split, bipartite and line graphs~\cite{ELW15}, unit disk graphs~\cite{gabriel}, interval and permutation graphs of diameter~2~\cite{FMNPV17}, and graphs of pathwidth~24~\cite{LP22}. On the positive side, it can easily be solved in linear time on trees~\cite{CEJO00,HM76,KRR96,S75}.
More involved polynomial-time algorithms exist for unicyclic graphs~\cite{SedlarUnicyclic} and, more generally, graphs of bounded cyclomatic number~\cite{ELW15}, outerplanar graphs~\cite{MDplanar}, cographs~\cite{ELW15}, chain graphs~\cite{FHHMS15}, cactus-block graphs~\cite{cactus}, and bipartite distance-hereditary graphs~\cite{bdh}. There are fixed parameter tractable (FPT) algorithms for the undirected graph parameters max leaf number~\cite{E15}, tree-depth~\cite{TD}, modular-width~\cite{MDwidth} and distance to cluster~\cite{GKIST22}, but FPT algorithms are highly unlikely to exist for the parameters solution size~\cite{HN13} and feedback vertex set~\cite{GKIST22}.

Due to the interest for \MD on undirected graphs, it is natural to ask what can be said in the context of digraphs. The metric dimension of digraphs was first studied in~\cite{CRZ00} under a somewhat restrictive definition; for our definitions, we follow the recent paper~\cite{ABCHONS20}, in which the algorithmic aspects of \MD on digraphs have been addressed. We call \emph{oriented graph} a digraph without directed 2-cycles. A \emph{directed acyclic digraph} (DAG for short) has no directed cycles at all. The \emph{underlying multigraph} of a digraph is the one obtained by ignoring the arc orientations; its \emph{underlying graph} is obtained from it by ignoring multiple edges. In a digraph, a \emph{strongly connected component} is a subgraph where every vertex is reachable from all other vertices. 
Note that for the \MD problem, undirected graphs can be seen as a special type of digraphs where each arc has a symmetric arc (\emph{i.e.}, replace every edge of the undirected graph by a directed 2-cycle).

The NP-hardness of \MD was proven for oriented graphs in~\cite{RRCM14} and, more recently, for bipartite DAGs of maximum degree~8 and maximum distance~4~\cite{ABCHONS20} (the \emph{maximum distance} being the length of a longest directed path without shortcuts). A linear-time algorithm for \MD on oriented trees was given in~\cite{ABCHONS20}.

\smallskip

\noindent\emph{Our results.} We generalize the linear-time algorithm for \MD on oriented trees from~\cite{ABCHONS20} to all digraphs whose underlying graph is a tree. In other words, here we allow 2-cycles. This makes a significant difference with oriented trees, and as a result our algorithm is non-trivial. We then extend the used methods to solve \MD in linear time for unicyclic digraphs (digraphs with a unique cycle). Then, we prove that \MD can be solved in time $f(t)n^{O(1)}$ for digraphs of order $n$ and modular-width~$t$ (a parameter recently introduced for digraphs in~\cite{SW20}). This extends the same result for undirected graphs from~\cite{MDwidth}, and is the first FPT algorithm for \MD on digraphs. Finally, we complement the hardness result from~\cite{ABCHONS20} by showing that \MD is NP-hard even for planar triangle-free DAGs of maximum degree 6 and maximum distance~4. 

A short preliminary version of this paper has appeared in the proceedings of the WG 2023 conference~\cite{WG}.

\section{Digraphs whose underlying graph is a tree (di-trees)}
\label{section-MDTree}

For the sake of convenience, we call \emph{di-tree} a digraph whose underlying graph is a tree. Trees are often the first non-trivial class to study for a graph problem. \MD{} is no exception to this, having been studied in the first papers for the undirected~\cite{CEJO00,HM76,KRR96,S75} and the oriented~\cite{ABCHONS20} cases. In the undirected case, a minimum-size resolving set can be found by taking, for each vertex of degree at least~3 spanning $k$ legs, the endpoint of $k-1$ of its legs (a \emph{leg} is an induced path spanning from a vertex of degree at least~3, having its inner vertices of degree~2, and ending in a leaf). In the case of oriented trees, taking all the sources (a \emph{source} is a vertex with no in-neighbor) and $k-1$ vertices in each set of $k$ in-twins yields a metric basis (two vertices are \emph{in-twins} if they have the same in-neighborhood). Our algorithm, being on di-trees (which include both undirected trees and oriented trees), will reuse those strategies, but we will need to refine them in order to obtain a metric basis.
The first refinement is of the notion of in-twins, for which we need the following notion:

\begin{definition}
	\label{def-escalator}
	A strongly connected component $E$ of a di-tree is an \emph{escalator} if it satisfies the following conditions:
	\begin{enumerate}
		\item its underlying graph is a path with vertices $e_1,\ldots,e_k$ ($k \geq 2$);
		\item there is a unique vertex $y \not\in E$ such that the arc $\overrightarrow{ye_1}$ exists and for all $i\in \{2,\ldots,k\}$, $e_i$ has no in-neighbors from outside E;
		\item for every $i \in \{1,\ldots,k-1\}$, no arc $\overrightarrow{e_iz}$ with $z \not\in E$ exists.
	\end{enumerate}
\end{definition}

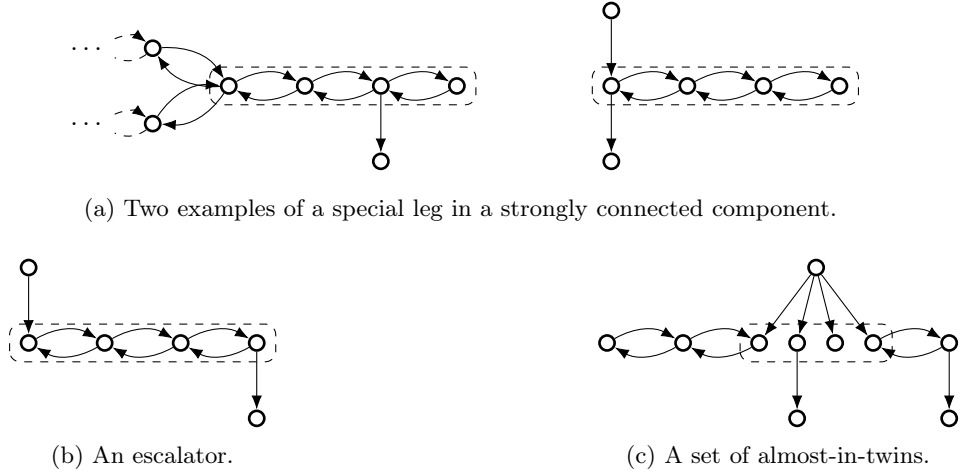
\begin{figure}[!ht]
	\centering
	\begin{subfigure}[b]{\textwidth}
		\centering
		\begin{tikzpicture}
			\node (a) at (0,0) {
				\begin{tikzpicture}
					\node[noeud] (ua) at (-1,0.5) {};
					\node[noeud] (ub) at (-1,-0.5) {};
					\node[noeud] (u0) at (0,0) {};
					\node[noeud] (u1) at (1,0) {};
					\node[noeud] (u2) at (2,0) {};
					\node[noeud] (u3) at (3,0) {};
					\draw[\nicearrow,bend left,dashed] (-1.5,0.625)to(ua);
					\draw[bend left,dashed] (ua)to(-1.5,0.375);
					\draw (-1.875,0.5) node {$\ldots$};
					\draw[\nicearrow,bend left,dashed] (-1.5,-0.375)to(ub);
					\draw[bend left,dashed] (ub)to(-1.5,-0.625);
					\draw (-1.875,-0.5) node {$\ldots$};
					\draw[\nicearrow,bend left] (u0)to(ua);
					\draw[\nicearrow,bend left] (ua)to(u0);
					\draw[\nicearrow,bend left] (u0)to(ub);
					\draw[\nicearrow,bend left] (ub)to(u0);
					\draw[\nicearrow,bend left] (u0)to(u1);
					\draw[\nicearrow,bend left] (u1)to(u0);
					\draw[\nicearrow,bend left] (u2)to(u1);
					\draw[\nicearrow,bend left] (u1)to(u2);
					\draw[\nicearrow,bend left] (u2)to(u3);
					\draw[\nicearrow,bend left] (u3)to(u2);
					\node (v0) at (0,1) {};
					\node[noeud] (v2) at (2,-1) {};
					\draw[\nicearrow] (u2)to(v2);
					\draw[dashed,rounded corners] (-0.25,-0.25) rectangle (3.25,0.25);
					\draw (u0) node[above,yshift=2.5mm] {$v$};
					\draw (u3) node[above,yshift=2.5mm] {$w$};
					\draw (u2) node[above,yshift=2.5mm] {$y$};
					\draw (v2) node[right,xshift=1mm] {$z$};
				\end{tikzpicture}
			};
			
			\node (b) at (6,0) {
				\begin{tikzpicture}
					\node[noeud] (u0) at (0,0) {};
					\node[noeud] (u1) at (1,0) {};
					\node[noeud] (u2) at (2,0) {};
					\node[noeud] (u3) at (3,0) {};
					\draw[\nicearrow,bend left] (u0)to(u1);
					\draw[\nicearrow,bend left] (u1)to(u0);
					\draw[\nicearrow,bend left] (u2)to(u1);
					\draw[\nicearrow,bend left] (u1)to(u2);
					\draw[\nicearrow,bend left] (u2)to(u3);
					\draw[\nicearrow,bend left] (u3)to(u2);
					\node[noeud] (v0) at (0,1) {};
					\node[noeud] (v0b) at (0,-1) {};
					\draw[\nicearrow] (v0)to(u0);
					\draw[\nicearrow] (u0)to(v0b);
					\draw[dashed,rounded corners] (-0.25,-0.25) rectangle (3.25,0.25);
					\draw (u0) node[left,xshift=-3mm] {$v=y$};
					\draw (u3) node[above,yshift=2.5mm] {$w$};
					\draw (v0b) node[right,xshift=1mm] {$z$};
				\end{tikzpicture}
			};
		\end{tikzpicture}
		\caption{Two examples of a special leg in a strongly connected component.}
		\label{fig-specialLeg}
	\end{subfigure}
	\vskip\baselineskip
	\begin{subfigure}[b]{0.475\textwidth}
		\centering
		\begin{tikzpicture}
			\node[noeud] (u0) at (0,0) {};
			\node[noeud] (u1) at (1,0) {};
			\node[noeud] (u2) at (2,0) {};
			\node[noeud] (u3) at (3,0) {};
			\draw[\nicearrow,bend left] (u0)to(u1);
			\draw[\nicearrow,bend left] (u1)to(u0);
			\draw[\nicearrow,bend left] (u2)to(u1);
			\draw[\nicearrow,bend left] (u1)to(u2);
			\draw[\nicearrow,bend left] (u2)to(u3);
			\draw[\nicearrow,bend left] (u3)to(u2);
			\node[noeud] (v0) at (0,1) {};
			\node[noeud] (v3) at (3,-1) {};
			\draw[\nicearrow] (v0)to(u0);
			\draw[\nicearrow] (u3)to(v3);
			\draw[dashed,rounded corners] (-0.25,-0.25) rectangle (3.25,0.25);
			\draw (u0) node[left,xshift=-2.5mm] {$e_1$};
			\draw (u1) node[above,yshift=2.5mm] {$e_2$};
			\draw (u2) node[above,yshift=2.5mm] {$e_3$};
			\draw (u3) node[right,xshift=2.5mm] {$e_4$};
			\draw (v0) node[left,xshift=-1mm] {$y$};
			\draw (v3) node[right,xshift=1mm] {$z$};
		\end{tikzpicture}
		\caption{An escalator.}
		\label{fig-escalator}
	\end{subfigure}
	\hfill
	\begin{subfigure}[b]{0.475\textwidth}
		\centering
		\begin{tikzpicture}
			\node[noeud] (u0) at (0,0) {};
			\node[noeud] (u1) at (1,0) {};
			\node[noeud] (u2) at (2,0) {};
			\node[noeud] (u25) at (3,0) {};
			\node[noeud] (u3) at (4,0) {};
			\node[noeud] (u35) at (5,0) {};
			\node[noeud] (u45) at (6,0) {};
			\draw[\nicearrow,bend left] (u0)to(u1);
			\draw[\nicearrow,bend left] (u1)to(u0);
			\draw[\nicearrow,bend left] (u2)to(u1);
			\draw[\nicearrow,bend left] (u1)to(u2);
			\draw[\nicearrow,bend left] (u45)to(u35);
			\draw[\nicearrow,bend left] (u35)to(u45);
			\node[noeud] (v275) at (3.5,1) {};
			\node[noeud] (v25) at (3,-1) {};
			\node[noeud] (v45) at (6,-1) {};
			\draw[\nicearrow] (v275)to(u2);
			\draw[\nicearrow] (v275)to(u25);
			\draw[\nicearrow] (v275)to(u3);
			\draw[\nicearrow] (v275)to(u35);
			\draw[\nicearrow] (u25)to(v25);
			\draw[\nicearrow] (u45)to(v45);
			\draw[dashed,rounded corners] (1.75,-0.25) rectangle (5.25,0.25);
			\draw (u2) node[below,yshift=-2.5mm] {$a_1$};
			\draw (u1) node[below,yshift=-2.5mm] {$y$};
			\draw (u25) node[left,xshift=-0.5mm] {$a_2$};
			\draw (u3) node[below,yshift=-2.5mm] {$a_3$};
			\draw (u35) node[below,yshift=-2.5mm] {$a_4$};
			\draw (u45) node[right,xshift=2.5mm] {$y$};
			\draw (v275) node[above,yshift=1.5mm] {$x$};
		\end{tikzpicture}
		\caption{A set of almost-in-twins.}
		\label{fig-almostInTwins}
	\end{subfigure}
	\caption{Illustrations of \Cref{def-specialLeg,def-escalator,def-almostInTwins}. The vertex names are taken from those definitions.}
	\label{fig-definitionsForTrees}
\end{figure}

An example of an escalator is depicted on \Cref{fig-escalator}.
Note that there can be any number (possibly, zero) of vertices $z \not\in E$ such that the arc $\overrightarrow{e_kz}$ exists.

\begin{definition}
	\label{def-almostInTwins}
	In a di-tree, a set of vertices $A=\{a_1,\ldots,a_k\}$ is a set of \emph{almost-in-twins} if there is a vertex $x$ such that:
	\begin{enumerate}
		\item for every $i \in \{1,\ldots,k\}$, the arc $\overrightarrow{xa_i}$ exists and the arc $\overrightarrow{a_ix}$ does not exist;
		\item for every $i \in \{1,\ldots,k\}$, either $a_i$ is a trivial strongly connected component and $N^-(a_i)=\{x\}$, or $a_i$ is the endpoint of an escalator and $N^-(a_i)=\{x,y\}$ where $y$ is its neighbor in the escalator.
	\end{enumerate}
\end{definition}

An example of a set of almost-in-twins is depicted on \Cref{fig-almostInTwins}.
Note that regular in-twins are also almost-in-twins. The second refinement is the following (for a given vertex $x$ in a strongly connected component with $C$ as an underlying graph, we call $d_C(x)$ the degree of $x$ in $C$):

\begin{definition}
	\label{def-specialLeg}
	Given the underlying graph $C$ of a strongly connected component of a di-tree and a set $D$ of vertices, we call a set $S$ of vertices inducing a path of order at least~2 in $C$ a \emph{special leg} if it verifies the four following properties:
	\begin{enumerate}
		\item $S$ has a unique vertex $v$ such that $v \in D$ or $d_C(v) \geq 3$;
		\item $S$ has a unique vertex $w$ such that $d_C(w) = 1$, furthermore $w \not\in D$: $w$ is called the \emph{endpoint} of $S$;
		\item all of the other vertices $x$ of $S$ verify $d_C(x) = 2$ and $x \not\in D$;
		\item at least one of the vertices $y \in S \setminus \{w\}$ has an out-arc $\overrightarrow{yz}$ with $z \notin C$ and $N^-(z)=\{y\}$.
	\end{enumerate}
\end{definition}

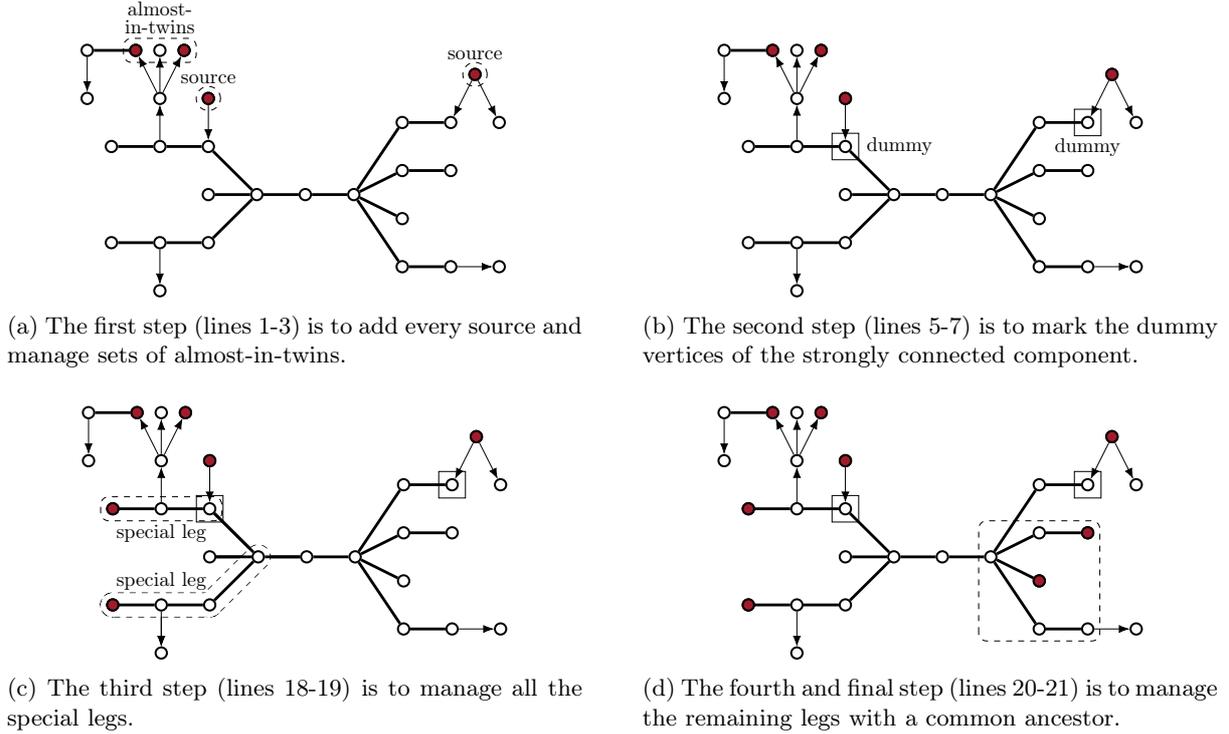
\begin{figure}[!ht]
	\centering
	
	\begin{subfigure}[b]{0.475\textwidth}
		\centering
		\scalebox{0.75}{
			\begin{tikzpicture}[scale=0.85]
				\node[noeud] (0) at (3,2) {};
				\node[noeud] (1) at (4,2) {};
				\node[noeud] (2) at (5,2) {};
				\draw[line width=0.5mm] (0)to(1)to(2);
				\node[noeud] (ul0) at (2,3) {};
				\node[noeud] (ul1) at (1,3) {};
				\node[noeud] (ul2) at (0,3) {};
				\node[noeud] (ul3) at (1,4) {};
				\node[noeud] (ul4) at (0.5,5) {};
				\node[noeud] (ul5) at (1,5) {};
				\node[noeud] (ul6) at (1.5,5) {};
				\node[noeud] (ul7) at (2,4) {};
				\node[noeud] (ul8) at (-0.5,5) {};
				\node[noeud] (ul9) at (-0.5,4) {};
				\draw[line width=0.5mm] (ul4)to(ul8);
				\draw[\nicearrow] (ul8)to(ul9);
				\draw[line width=0.5mm] (0)to(ul0)to(ul1)to(ul2);
				\draw[\nicearrow] (ul1)to(ul3);
				\draw[\nicearrow] (ul3)to(ul4);
				\draw[\nicearrow] (ul3)to(ul5);
				\draw[\nicearrow] (ul3)to(ul6);
				\draw[\nicearrow] (ul7)to(ul0);
				\node[noeud] (ml0) at (2,2) {};
				\draw[line width=0.5mm] (0)to(ml0);
				\node[noeud] (ll0) at (2,1) {};
				\node[noeud] (ll1) at (1,1) {};
				\node[noeud] (ll2) at (0,1) {};
				\node[noeud] (ll3) at (1,0) {};
				\draw[line width=0.5mm] (0)to(ll0)to(ll1)to(ll2);
				\draw[\nicearrow] (ll1)to(ll3);
				\node[noeud] (ur0) at (6,3.5) {};
				\node[noeud] (ur1) at (7,3.5) {};
				\node[noeud] (ur2) at (7.5,4.5) {};
				\node[noeud] (ur3) at (8,3.5) {};
				\draw[line width=0.5mm] (2)to(ur0)to(ur1);
				\draw[\nicearrow] (ur2)to(ur1);
				\draw[\nicearrow] (ur2)to(ur3);
				\node[noeud] (mr0) at (6,2.5) {};
				\node[noeud] (mr1) at (7,2.5) {};
				\node[noeud] (mr2) at (6,1.5) {};
				\node[noeud] (mr3) at (6,0.5) {};
				\node[noeud] (mr4) at (7,0.5) {};
				\node[noeud] (mr5) at (8,0.5) {};
				\draw[line width=0.5mm] (2)to(mr0)to(mr1);
				\draw[line width=0.5mm] (2)to(mr2);
				\draw[line width=0.5mm] (2)to(mr3)to(mr4);
				\draw[\nicearrow] (mr4)to(mr5);
				
				\node[noeud,fill=nicered] (ur2) at (7.5,4.5) {};
				\node[noeud,fill=nicered] (ul6) at (2,4) {};
				\draw[dashed] (ur2) circle (0.25);
				\node at (ur2) [above=1.5mm] {source};
				\draw[dashed] (ul6) circle (0.25);
				\node at (ul6) [above=1.5mm] {source};
				
				\node[noeud,fill=nicered] (ul4) at (0.5,5) {};
				\node[noeud,fill=nicered] (ul6) at (1.5,5) {};
				\draw[dashed,rounded corners] (0.25,4.75) rectangle (1.75,5.25);
				\node at (1,5.875) {almost-};
				\node at (1,5.5) {in-twins};
			\end{tikzpicture}
		}
		\caption{The first step (lines 1-3) is to add every source and manage sets of almost-in-twins.}
		\label{fig-MDTree-SourceInTwins}
	\end{subfigure}
	\hfill
	\begin{subfigure}[b]{0.475\textwidth}
		\centering
		\scalebox{0.75}{
			\begin{tikzpicture}[scale=0.85]
				\node[noeud] (0) at (3,2) {};
				\node[noeud] (1) at (4,2) {};
				\node[noeud] (2) at (5,2) {};
				\draw[line width=0.5mm] (0)to(1)to(2);
				\node[noeud] (ul0) at (2,3) {};
				\node[noeud] (ul1) at (1,3) {};
				\node[noeud] (ul2) at (0,3) {};
				\node[noeud] (ul3) at (1,4) {};
				\node[noeud] (ul4) at (0.5,5) {};
				\node[noeud] (ul5) at (1,5) {};
				\node[noeud] (ul6) at (1.5,5) {};
				\node[noeud] (ul7) at (2,4) {};
				\node[noeud] (ul8) at (-0.5,5) {};
				\node[noeud] (ul9) at (-0.5,4) {};
				\draw[line width=0.5mm] (ul4)to(ul8);
				\draw[\nicearrow] (ul8)to(ul9);
				\draw[line width=0.5mm] (0)to(ul0)to(ul1)to(ul2);
				\draw[\nicearrow] (ul1)to(ul3);
				\draw[\nicearrow] (ul3)to(ul4);
				\draw[\nicearrow] (ul3)to(ul5);
				\draw[\nicearrow] (ul3)to(ul6);
				\draw[\nicearrow] (ul7)to(ul0);
				\node[noeud] (ml0) at (2,2) {};
				\draw[line width=0.5mm] (0)to(ml0);
				\node[noeud] (ll0) at (2,1) {};
				\node[noeud] (ll1) at (1,1) {};
				\node[noeud] (ll2) at (0,1) {};
				\node[noeud] (ll3) at (1,0) {};
				\draw[line width=0.5mm] (0)to(ll0)to(ll1)to(ll2);
				\draw[\nicearrow] (ll1)to(ll3);
				\node[noeud] (ur0) at (6,3.5) {};
				\node[noeud] (ur1) at (7,3.5) {};
				\node[noeud] (ur2) at (7.5,4.5) {};
				\node[noeud] (ur3) at (8,3.5) {};
				\draw[line width=0.5mm] (2)to(ur0)to(ur1);
				\draw[\nicearrow] (ur2)to(ur1);
				\draw[\nicearrow] (ur2)to(ur3);
				\node[noeud] (mr0) at (6,2.5) {};
				\node[noeud] (mr1) at (7,2.5) {};
				\node[noeud] (mr2) at (6,1.5) {};
				\node[noeud] (mr3) at (6,0.5) {};
				\node[noeud] (mr4) at (7,0.5) {};
				\node[noeud] (mr5) at (8,0.5) {};
				\draw[line width=0.5mm] (2)to(mr0)to(mr1);
				\draw[line width=0.5mm] (2)to(mr2);
				\draw[line width=0.5mm] (2)to(mr3)to(mr4);
				\draw[\nicearrow] (mr4)to(mr5);
				
				\node[noeud,fill=nicered] (ur2) at (7.5,4.5) {};
				\node[noeud,fill=nicered] (ul6) at (2,4) {};
				\node[noeud,fill=nicered] (ul4) at (0.5,5) {};
				\node[noeud,fill=nicered] (ul6) at (1.5,5) {};
				
				\node[draw,rectangle,fit=(ul0)] {};
				\node at (ul0) [right=2.5mm] {dummy};
				\node[draw,rectangle,fit=(ur1)] {};
				\node at (ur1) [below=1.5mm] {dummy};
			\end{tikzpicture}
		}
		\caption{The second step (lines 5-7) is to mark the dummy vertices of the strongly connected component.}
		\label{fig-MDTree-DummyVertices}
	\end{subfigure}
	\vskip\baselineskip
	\begin{subfigure}[b]{0.475\textwidth}
		\centering
		\scalebox{0.75}{
			\begin{tikzpicture}[scale=0.85]
				\node[noeud] (0) at (3,2) {};
				\node[noeud] (1) at (4,2) {};
				\node[noeud] (2) at (5,2) {};
				\draw[line width=0.5mm] (0)to(1)to(2);
				\node[noeud] (ul0) at (2,3) {};
				\node[noeud] (ul1) at (1,3) {};
				\node[noeud] (ul2) at (0,3) {};
				\node[noeud] (ul3) at (1,4) {};
				\node[noeud] (ul4) at (0.5,5) {};
				\node[noeud] (ul5) at (1,5) {};
				\node[noeud] (ul6) at (1.5,5) {};
				\node[noeud] (ul7) at (2,4) {};
				\node[noeud] (ul8) at (-0.5,5) {};
				\node[noeud] (ul9) at (-0.5,4) {};
				\draw[line width=0.5mm] (ul4)to(ul8);
				\draw[\nicearrow] (ul8)to(ul9);
				\draw[line width=0.5mm] (0)to(ul0)to(ul1)to(ul2);
				\draw[\nicearrow] (ul1)to(ul3);
				\draw[\nicearrow] (ul3)to(ul4);
				\draw[\nicearrow] (ul3)to(ul5);
				\draw[\nicearrow] (ul3)to(ul6);
				\draw[\nicearrow] (ul7)to(ul0);
				\node[noeud] (ml0) at (2,2) {};
				\draw[line width=0.5mm] (0)to(ml0);
				\node[noeud] (ll0) at (2,1) {};
				\node[noeud] (ll1) at (1,1) {};
				\node[noeud] (ll2) at (0,1) {};
				\node[noeud] (ll3) at (1,0) {};
				\draw[line width=0.5mm] (0)to(ll0)to(ll1)to(ll2);
				\draw[\nicearrow] (ll1)to(ll3);
				\node[noeud] (ur0) at (6,3.5) {};
				\node[noeud] (ur1) at (7,3.5) {};
				\node[noeud] (ur2) at (7.5,4.5) {};
				\node[noeud] (ur3) at (8,3.5) {};
				\draw[line width=0.5mm] (2)to(ur0)to(ur1);
				\draw[\nicearrow] (ur2)to(ur1);
				\draw[\nicearrow] (ur2)to(ur3);
				\node[noeud] (mr0) at (6,2.5) {};
				\node[noeud] (mr1) at (7,2.5) {};
				\node[noeud] (mr2) at (6,1.5) {};
				\node[noeud] (mr3) at (6,0.5) {};
				\node[noeud] (mr4) at (7,0.5) {};
				\node[noeud] (mr5) at (8,0.5) {};
				\draw[line width=0.5mm] (2)to(mr0)to(mr1);
				\draw[line width=0.5mm] (2)to(mr2);
				\draw[line width=0.5mm] (2)to(mr3)to(mr4);
				\draw[\nicearrow] (mr4)to(mr5);
				
				\node[noeud,fill=nicered] (ur2) at (7.5,4.5) {};
				\node[noeud,fill=nicered] (ul6) at (2,4) {};
				\node[noeud,fill=nicered] (ul4) at (0.5,5) {};
				\node[noeud,fill=nicered] (ul6) at (1.5,5) {};
				
				\node[draw,rectangle,fit=(ul0)] {};
				\node[draw,rectangle,fit=(ur1)] {};
				
				\draw[dashed,rounded corners] (-0.25,2.75) rectangle (2.25,3.25);
				\node at (1,2.5) {special leg};
				\draw[dashed] \convexpath{ll2,ll1,ll0,0,ll0,ll1,ll2}{0.25cm};
				\fill[white] \convexpath{ll2,ll1,ll0,0,ll0,ll1,ll2}{0.2425cm};
				\node at (1,1.5) {special leg};
				\node[noeud,fill=nicered] (ul2) at (0,3) {};
				\node[noeud,fill=nicered] (ll2) at (0,1) {};
				
				\node[noeud] (0) at (3,2) {};
				\node[noeud] (ll0) at (2,1) {};
				\node[noeud] (ll1) at (1,1) {};
				\node[noeud,fill=nicered] (ll2) at (0,1) {};
				\draw[line width=0.5mm] (0)to(1);
				\draw[line width=0.5mm] (0)to(ul0);
				\draw[line width=0.5mm] (0)to(ml0);
				\draw[line width=0.5mm] (0)to(ll0)to(ll1)to(ll2);
				\draw[\nicearrow] (ll1)to(ll3);
			\end{tikzpicture}
		}
		\caption{The third step (lines 18-19) is to manage all the special legs.}
		\label{fig-MDTree-SpecialLegs}
	\end{subfigure}
	\hfill
	\begin{subfigure}[b]{0.475\textwidth}
		\centering
		\scalebox{0.75}{
			\begin{tikzpicture}[scale=0.85]
				\node[noeud] (0) at (3,2) {};
				\node[noeud] (1) at (4,2) {};
				\node[noeud] (2) at (5,2) {};
				\draw[line width=0.5mm] (0)to(1)to(2);
				\node[noeud] (ul0) at (2,3) {};
				\node[noeud] (ul1) at (1,3) {};
				\node[noeud] (ul2) at (0,3) {};
				\node[noeud] (ul3) at (1,4) {};
				\node[noeud] (ul4) at (0.5,5) {};
				\node[noeud] (ul5) at (1,5) {};
				\node[noeud] (ul6) at (1.5,5) {};
				\node[noeud] (ul7) at (2,4) {};
				\node[noeud] (ul8) at (-0.5,5) {};
				\node[noeud] (ul9) at (-0.5,4) {};
				\draw[line width=0.5mm] (ul4)to(ul8);
				\draw[\nicearrow] (ul8)to(ul9);
				\draw[line width=0.5mm] (0)to(ul0)to(ul1)to(ul2);
				\draw[\nicearrow] (ul1)to(ul3);
				\draw[\nicearrow] (ul3)to(ul4);
				\draw[\nicearrow] (ul3)to(ul5);
				\draw[\nicearrow] (ul3)to(ul6);
				\draw[\nicearrow] (ul7)to(ul0);
				\node[noeud] (ml0) at (2,2) {};
				\draw[line width=0.5mm] (0)to(ml0);
				\node[noeud] (ll0) at (2,1) {};
				\node[noeud] (ll1) at (1,1) {};
				\node[noeud] (ll2) at (0,1) {};
				\node[noeud] (ll3) at (1,0) {};
				\draw[line width=0.5mm] (0)to(ll0)to(ll1)to(ll2);
				\draw[\nicearrow] (ll1)to(ll3);
				\node[noeud] (ur0) at (6,3.5) {};
				\node[noeud] (ur1) at (7,3.5) {};
				\node[noeud] (ur2) at (7.5,4.5) {};
				\node[noeud] (ur3) at (8,3.5) {};
				\draw[line width=0.5mm] (2)to(ur0)to(ur1);
				\draw[\nicearrow] (ur2)to(ur1);
				\draw[\nicearrow] (ur2)to(ur3);
				\node[noeud] (mr0) at (6,2.5) {};
				\node[noeud] (mr1) at (7,2.5) {};
				\node[noeud] (mr2) at (6,1.5) {};
				\node[noeud] (mr3) at (6,0.5) {};
				\node[noeud] (mr4) at (7,0.5) {};
				\node[noeud] (mr5) at (8,0.5) {};
				\draw[line width=0.5mm] (2)to(mr0)to(mr1);
				\draw[line width=0.5mm] (2)to(mr2);
				\draw[line width=0.5mm] (2)to(mr3)to(mr4);
				\draw[\nicearrow] (mr4)to(mr5);
				
				\node[noeud,fill=nicered] (ur2) at (7.5,4.5) {};
				\node[noeud,fill=nicered] (ul6) at (2,4) {};
				\node[noeud,fill=nicered] (ul4) at (0.5,5) {};
				\node[noeud,fill=nicered] (ul6) at (1.5,5) {};
				\node[noeud,fill=nicered] (ul2) at (0,3) {};
				\node[noeud,fill=nicered] (ll2) at (0,1) {};
				
				\node[draw,rectangle,fit=(ul0)] {};
				\node[draw,rectangle,fit=(ur1)] {};
				
				\draw[dashed,rounded corners] (4.75,0.25) rectangle (7.25,2.75);
				\node[noeud,fill=nicered] (mr1) at (7,2.5) {};
				\node[noeud,fill=nicered] (mr2) at (6,1.5) {};
			\end{tikzpicture}
		}
		\caption{The fourth and final step (lines 20-21) is to manage the remaining legs with a common ancestor.}
		\label{fig-MDTree-RemainingLegs}
	\end{subfigure}
	
	\caption{Illustration of \Cref{alg-MDTree}. For the sake of simplicity, there are only two strongly connected components, for which we only represent the underlying graph with bolded edges, so every bolded edge is a 2-cycle. One of the two strongly connected components is a simple path that does not require any action. Vertices in the metric basis are colored in red.}
	\label{fig-MDTree}
\end{figure}

\begin{figure}[!ht]
	\centering
	\begin{subfigure}[b]{0.475\textwidth}
		\centering
		\begin{tikzpicture}
			\node[noeud,fill=nicered] (u) at (0,0) {};
			\node[noeud] (v) at (2,0) {};
			\draw[line width=0.5mm,style={decorate, decoration=snake}] (u)to(v);
			\node (a) at (1,-1) {};
		\end{tikzpicture}
		\caption{Neither in-arc nor out-arc.}
	\end{subfigure}
	\hfill
	\begin{subfigure}[b]{0.475\textwidth}
		\centering
		\begin{tikzpicture}
			\node (p1) at (0,1.5) {
				\begin{tikzpicture}
					\node[noeud] (u) at (0,0) {};
					\node[noeud,fill=nicered] (v) at (2,0) {};
					\node[noeud] (w) at (0,-0.75) {};
					\draw[\nicearrow] (u)to(w);
					\draw[line width=0.5mm,style={decorate, decoration=snake}] (u)to(v);
				\end{tikzpicture}
			};
			\node (p2) at (0,0) {
				\begin{tikzpicture}
					\node[noeud,fill=nicered] (u) at (0,0) {};
					\node[noeud] (v) at (2,0) {};
					\node[noeud] (w) at (2,-0.75) {};
					\draw[\nicearrow] (v)to(w);
					\draw[line width=0.5mm,style={decorate, decoration=snake}] (u)to(v);
				\end{tikzpicture}
			};
		\end{tikzpicture}
		\caption{No in-arc, out-arc only at an endpoint.}
	\end{subfigure}
	\vskip\baselineskip
	\begin{subfigure}[b]{0.475\textwidth}
		\centering
		\begin{tikzpicture}
			\node (p1) at (0,1.5) {
				\begin{tikzpicture}
					\node[noeud,fill=nicered] (u) at (0,0) {};
					\node[noeud,fill=nicered] (v) at (2,0) {};
					\node[noeud] (w) at (1,0) {};
					\node[noeud] (wp) at (1,-0.75) {};
					\draw[\nicearrow] (w)to(wp);
					\draw[line width=0.5mm,style={decorate, decoration=snake}] (u)to(w)to(v);
				\end{tikzpicture}
			};
			\node (p2) at (0,0) {
				\begin{tikzpicture}
					\node[noeud,fill=nicered] (u) at (0,0) {};
					\node[noeud,fill=nicered] (v) at (2,0) {};
					\node[noeud] (w) at (1,0) {};
					\node[noeud] (wp) at (1,-0.75) {};
					\node[noeud] (up) at (0,-0.75) {};
					\draw[\nicearrow] (w)to(wp);
					\draw[\nicearrow] (u)to(up);
					\draw[line width=0.5mm,style={decorate, decoration=snake}] (u)to(w)to(v);
				\end{tikzpicture}
			};
			\node (p3) at (2.5,1.5) {
				\begin{tikzpicture}
					\node[noeud,fill=nicered] (u) at (0,0) {};
					\node[noeud,fill=nicered] (v) at (2,0) {};
					\node[noeud] (up) at (0,-0.75) {};
					\node[noeud] (vp) at (2,-0.75) {};
					\draw[\nicearrow] (u)to(up);
					\draw[\nicearrow] (v)to(vp);
					\draw[line width=0.5mm,style={decorate, decoration=snake}] (u)to(v);
				\end{tikzpicture}
			};
			\node (p4) at (2.5,0) {
				\begin{tikzpicture}
					\node[noeud,fill=nicered] (u) at (0,0) {};
					\node[noeud,fill=nicered] (v) at (2,0) {};
					\node[noeud] (w) at (1,0) {};
					\node[noeud] (wp) at (1,-0.75) {};
					\node[noeud] (up) at (0,-0.75) {};
					\node[noeud] (vp) at (2,-0.75) {};
					\draw[\nicearrow] (w)to(wp);
					\draw[\nicearrow] (u)to(up);
					\draw[\nicearrow] (v)to(vp);
					\draw[line width=0.5mm,style={decorate, decoration=snake}] (u)to(w)to(v);
				\end{tikzpicture}
			};
		\end{tikzpicture}
		\caption{No in-arc, other cases.}
	\end{subfigure}
	\hfill
	\begin{subfigure}[b]{0.475\textwidth}
		\centering
		\begin{tikzpicture}
			\node[noeud,fill=nicered] (u) at (0,0) {};
			\node[noeud] (v) at (2,0) {};
			\node[noeud] (w) at (1,0) {};
			\node[draw,rectangle,fit=(w)] {};
			\node[noeud] (wp) at (1,0.75) {};
			\draw[\nicearrow] (wp)to(w);
			\draw[line width=0.5mm,style={decorate, decoration=snake}] (u)to(w)to(v);
			\node (a) at (1,-1) {};
		\end{tikzpicture}
		\caption{One in-arc, in the middle.}
	\end{subfigure}
	\caption{The cases of \Cref{alg-MDTree} where a strongly connected component is a path and we have to add specific vertices (colored in red) to the metric basis. The path is depicted as wavy bolded edges. Note that the considered out-arcs have to be towards a vertex with no other in-neighbor than the one in the path (out-arcs to vertices with other in-neighbors can still exist).}
	\label{fig-MDTree-Path}
\end{figure}
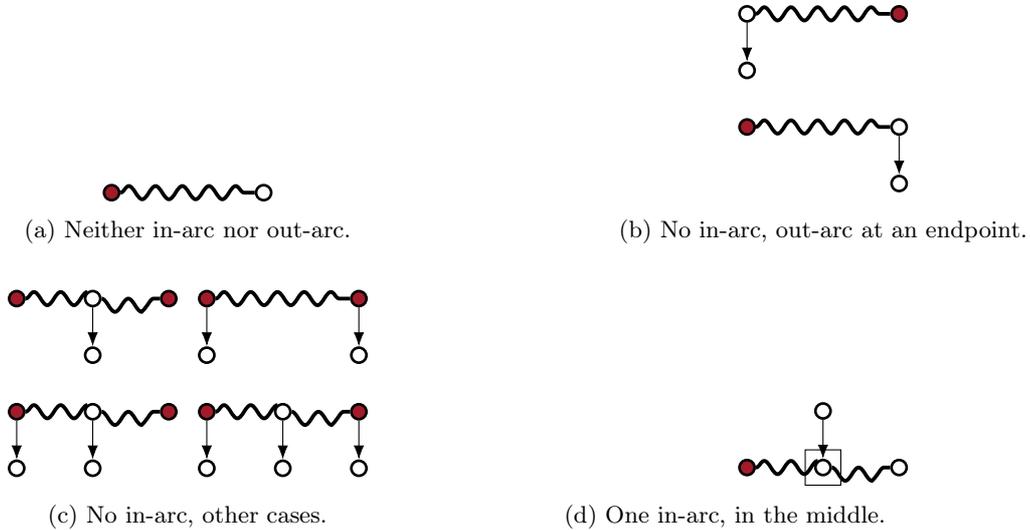

An example of a special leg is depicted on \Cref{fig-specialLeg}.
Note that several special legs can span from the same vertex, from which regular legs can also span.
\Cref{alg-MDTree}, illustrated in \Cref{fig-MDTree,fig-MDTree-Path}, computes a metric basis of a di-tree.

\begin{algorithm}[!ht]
	\caption{An algorithm computing the metric basis of a di-tree.}\label{alg-MDTree}
	\SetKwInOut{KwIn}{Input}
	\SetKwInOut{KwOut}{Output}
	\KwIn{A di-tree $T$.}
	\KwOut{A metric basis $\mathcal{B}$ of $T$.}
	
	$\mathcal{B}$ $\leftarrow$ Every source of $T$
	
	\ForEach{set $I$ of almost-in-twins}{
		Add $|I|-1$ vertices of $I$ to $\mathcal{B}$
	}
	
	\ForEach{strongly connected component with $C$ as an underlying graph}{
		$D$ $\leftarrow$ $\emptyset$\\
		\ForEach{arc $\overrightarrow{uv}$ with $v \in C$ and $u \not\in C$}{
			Add $v$ to $D$
		}
		
		\If{$C$ is a path with endpoints $x$ and $y$}{
			\If{there is no vertex in $C \cap D$}{
				\If{there is no arc $\overrightarrow{uv}$ such that $u \in C$, $v \not\in C$ and $N^-(v)=\{u\}$}{
					Add $x$ to $\mathcal{B}$
				}
				\ElseIf{the only arcs $\overrightarrow{uv}$ such that $u \in C$, $v \not\in C$ and $N^-(v)=\{u\}$ are such that $u=x$ (resp. $u=y$)}{
					Add $y$ (resp. $x$) to $\mathcal{B}$
				}
				\Else{
					Add $x$ and $y$ to $\mathcal{B}$
				}
			}
			\ElseIf{there is exactly one vertex $w$ in $C \cap D$, $w$ is neither $x$ nor $y$, and there is no arc $\overrightarrow{wz}$ such that $z \not\in C$ and $N^-(z)=\{w\}$}
			{
				Add $x$ to $\mathcal{B}$
			}
		}
		
		\ForEach{special leg $L$ of $C$}{
			Add the endpoint of $L$ to $\mathcal{B}$
		}
		
		\ForEach{vertex of degree $\geq 3$ in $C$ from which span $k \geq 2$ legs of $C$ that do not have a vertex in $\mathcal{B}$ or in $D$}{
			Add the endpoint of $k-1$ such legs to $\mathcal{B}$
		}
	}
	
	\Return $\mathcal{B}$
\end{algorithm}

\ \\\noindent\textbf{Explanation of \Cref{alg-MDTree}.} The algorithm will compute a metric basis $\mathcal{B}$ of a di-tree $T$ in linear time.
The first thing we do is to add every source in $T$ to $\mathcal{B}$ (line~1). Then, for every set of almost-in-twins, we add all of them but one to $\mathcal{B}$ (lines~2-3).
Those two first steps, depicted in \Cref{fig-MDTree-SourceInTwins}, are the ones used to compute the metric basis of an orientation of a tree~\cite{ABCHONS20}, and as such they are still necessary for managing the non-strongly connected components of the di-tree.
Note that we are specifically managing sets of \emph{almost-in-twins}, which include sets of in-twins, since it is necessary to resolve the specific case of escalators.
The rest of the algorithm consists in managing the strongly connected components.

For each strongly connected component having $C$ as an underlying graph, we first identify each vertex $x$ of $C$ that has an in-arc coming from \textbf{outside} $C$. Indeed, since $x$ is the "last" vertex of a path coming from outside $C$, there are vertices of $\mathcal{B}$ "behind" this in-arc (or they can themselves be a vertex in $\mathcal{B}$), which we will call $\mathcal{B}_x$. However, the vertices in $\mathcal{B}_x$ can be "projected" on $x$ since, $T$ being a di-tree, $x$ is on every shortest path from the vertices of $\mathcal{B}$ "behind" the in-arc to the vertices of $C$. Hence, we will mark $x$ as a \textbf{dummy vertex} (lines~5-7, depicted in \Cref{fig-MDTree-DummyVertices}): we will consider that it is in $\mathcal{B}$ for the rest of this step, and acts as a representative of the set $\mathcal{B}_x$ with respect to $C$.

We then have to manage some specific cases whenever $C$ is a path (lines 8-17). Indeed, the last two steps of the algorithm do not always work under some conditions. Those specific conditions will be highlighted in the proof, and are depicted in \Cref{fig-MDTree-Path}.

The last two steps are then applied. First, we have to consider the \textbf{special legs} defined in \Cref{def-specialLeg}. The idea behind those special legs is the following: for every out-arc $\overrightarrow{yz}$ with $y$ in the special leg and $z$ outside of $C$, any vertex in the metric basis "before" the start of the special leg will not distinguish $z$ and the next neighbor of $y$ in the special leg. Hence, we have to add at least one vertex to $\mathcal{B}$ for each special leg, and we choose the endpoint of the special leg (lines~18-19, depicted in \Cref{fig-MDTree-SpecialLegs}).
Finally, we apply the well-known algorithm for computing the metric basis of a tree to the remaining parts of~$C$ (lines~20-21, depicted in \Cref{fig-MDTree-RemainingLegs}). The special legs and the legs containing a dummy vertex, being already resolved, are not considered in this part.

\begin{theorem}
	\label{thm-MDTree}
	\Cref{alg-MDTree} computes a metric basis of a di-tree in linear time.
\end{theorem}

\begin{proof}
	Let $T$ be a di-tree, and $\mathcal{B}$ be the set of vertices returned by \Cref{alg-MDTree}. We need to prove that $\mathcal{B}$ resolves every pair of vertices of $T$, that $\mathcal{B}$ is of minimum size, and that \Cref{alg-MDTree} runs in linear time.
	
	First, note that, if $T$ is either an orientation of a tree or a strongly connected graph (and thus seen as an undirected graph), then, \Cref{alg-MDTree} does compute a metric basis. In the first case, $\mathcal{B}$ will contain only the vertices added in lines~1-3 (which correspond to a so-called \emph{adequate} set of vertices in~\cite{ABCHONS20}, see Lemma~2.10 and Theorem~2.11); and in the second case, $\mathcal{B}$ will contain either only a vertex added in lines~10-11 (if $T$ is a path) or only the vertices added in lines~20-21 (which correspond to the well-known resolving set of trees, see for example Theorem~2.4 in~\cite{KRR96}); those two cases do form a metric basis of $T$.
	
	Hence, we will consider that $T$ contains at least one strongly connected component and at least one non-strongly connected component.
	We will first show that the vertices we select in $\mathcal{B}$ are necessary to resolve at least one pair of vertices, and then that they do indeed form a resolving set, and thus that $\mathcal{B}$ is a metric basis.
	
	First, let us consider the sources. It is easy to see that each source has to be in $\mathcal{B}$, since otherwise they would not be reachable from any other vertex in $\mathcal{B}$.
	Now, let us consider the in-twins. Again, it is easy to see that the only way to resolve two in-twins will be to have at least one of them in $\mathcal{B}$. However, \Cref{alg-MDTree} considers \emph{almost-in-twins}, which are more general than regular in-twins. This is because of the \emph{escalators}: let $u$ be the endpoint of an escalator having an in-arc coming from outside of it. If $u$ has an almost-in-twin $v$, then, by definition, $u$ and $v$ cannot be resolved without taking in $\mathcal{B}$ either at least one of these or another vertex from the escalator. Thus, we choose to take either $u$ or $v$, which will resolve vertices in the escalator as well as the almost-in-twins.
	
	Note that sources and in-twins can only exist in a non-strongly connected component. Hence, in the rest of this part of the proof, we will consider a strongly connected component with $C$ as an underlying graph. Note that, by our construction, if a vertex $x \in C$ is a dummy vertex, then, there will be a vertex $b \in \mathcal{B}$ such that there is a path from $b$ to $x$ (even if it is not necessarily the case yet).
	
	We first consider the case of the special legs of $C$. Let $L$ be a special leg with vertices $x_1,\ldots,x_\ell$, starting from a vertex $x_1$ such that $d_C(x_1) \geq 3$ or $x_1$ is a dummy vertex, and ending at a vertex $x_\ell$ verifying $d_C(x_\ell)=1$. Let $x_i$ ($i \in \{1,\ldots,\ell\}$) be a vertex such that $d_C(x_i)\geq 2$ and there is an arc $\overrightarrow{x_iy}$ with $y \not\in C$ and $N^-(y)=\{x_i\}$. Note that, at this point, no vertex in the special leg can be in $\mathcal{B}$, and there is no arc from outside of $C$ to a vertex of $L$ apart from $x_1$. This means that, for vertices of $L$, $x_1$ can be seen as a representative from the set $\mathcal{B}$: every path from a vertex in $\mathcal{B}$ to a vertex in $L$ has to go through $x_1$. Furthermore, $y$, having no other in-neighbor, is only reached by vertices of $\mathcal{B}$ through $x_1$. In practice, this means that $x_{i+1}$ and $y$ are not resolved: since they are at the same distance from $x_1$, they are at the same distance from any vertex in $\mathcal{B}$. Hence, it is necessary to add either a vertex from the set $\{x_{i+1},\ldots,x_\ell\}$ or $y$ to $\mathcal{B}$ in order to resolve this pair. Since other such situations might occur in $L$, the easiest solution to resolve this pair is to add $x_\ell$ to $\mathcal{B}$, since doing so will resolve all such pairs.
	
	We then consider the rest of $C$. Since $T$ is a di-tree, $C$ is a tree. Let $L_1$ and $L_2$ be two (non-special) legs spanning from a vertex $x$, and assume that neither of those legs have an in-arc coming from outside of $C$. If no vertex from either $L_1$ or $L_2$ is in $\mathcal{B}$, then, the vertices from $L_1$ and $L_2$ cannot be resolved. Here, we add the endpoint of either $L_1$ or $L_2$ to $\mathcal{B}$. Note that, since $C$ is not a path (which either is a special case that we will consider below, or has been considered in the special leg case), if one of $L_1$ or $L_2$ has a dummy vertex or a vertex in $\mathcal{B}$, then, the vertices from $L_1$ and $L_2$ will be resolved without having to add another vertex in~$\mathcal{B}$.
	
	However, note that there are cases where the above method does not create a resolving set. Indeed, when $C$ is a path with vertices $x_1,\ldots,x_k$, it is possible to fall into some patterns where we need to add specific vertices to $\mathcal{B}$. The patterns are the following:
	\begin{enumerate}
		\item There is no arc $\overrightarrow{yx_i}$ with $y \not\in C$, in which case we have the following subpatterns depending on the presence of out-arcs $\overrightarrow{x_iy}$ with $y \not\in C$ and $N^-(y) = \{x_i\}$:
		\begin{enumerate}
			\item there is no such arc $\overrightarrow{x_iy}$, in which case we have to add either $x_1$ or $x_k$ to $\mathcal{B}$ (the vertices of $C$ have to be reached from $\mathcal{B}$ and separated);
			\item there is such an arc and we have $u=x_1$ (resp. $u=x_k$), in which case we have to add a vertex from $C$ to $\mathcal{B}$ in order to reach the vertices of $C$; we add $x_k$ (resp. $x_1$) to $\mathcal{B}$ since, otherwise, either $y$ and $x_2$ (resp. $x_{k-1}$) will not be resolved or some pair of vertices in the path will not be resolved;
			\item in every other case (that is, there are such arcs leaving $C$ from different vertices), we have to add both $x_1$ and $x_k$ to $\mathcal{B}$ (since, otherwise, one vertex from $C$ and one vertex reached by an out-arc from $C$ will not be resolved).
		\end{enumerate}
		Note that there can be arcs $\overrightarrow{x_iy}$ with $y \not\in C$ and where $y$ has an in-neighbor not in $C$; those do not matter for this specific analysis, since such $y$ are already resolved from the vertices in $C$ by vertices in $\mathcal{B}$ that are ancestors of their other in-neighbors.
		\item There is exactly one $x_i$, $i \not\in \{1,k\}$, such that there exists an arc $\overrightarrow{yx_i}$ with $y \not\in C$ and there is no arc $\overrightarrow{x_iz}$ with $z \not\in C$ and $N^-(z)=\{x_i\}$, in which case we add either endpoint to $\mathcal{B}$ (since, otherwise, $x_{i-1}$ and $x_{i+1}$ will not be resolved). Note that if such an arc $\overrightarrow{x_iz}$ exists, then we have two special legs spanning from $x_i$ and thus no problem arises.
	\end{enumerate}
	In every other case where $C$ is a path, the cases considered above (in particular, the special leg and the sources and almost-in-twins) will have us add to $\mathcal{B}$ the vertices necessary to resolve the vertices of $C$ and its direct out-neighborhood.
	
	\ \\
	
	At this point, every vertex that we added to $\mathcal{B}$ was necessary to resolve at least one pair of vertices which could not be resolved any other way. At each step, when we had the choice between several vertices, we chose as few as possible to resolve everything. Hence, we only need to check that every pair of vertices is resolved by $\mathcal{B}$, which will prove that $\mathcal{B}$ is a metric basis.
	
	Assume by contradiction that two vertices $u$ and $v$ are not resolved by $\mathcal{B}$. This means that, for every vertex $b \in \mathcal{B}$, either both $u$ and $v$ are not reachable from $b$, or there are two unique shortest paths $P^b_u$ between $b$ and $u$ and $P^b_v$ between $b$ and $v$ such that $|P^b_u|=|P^b_v|$. Note that $P^b_u$ and $P^b_v$ are unique because the underlying graph of $T$ is a tree. Given a vertex $b \in \mathcal{B}$, let $x$ be the last common vertex of $P^b_u$ and $P^b_v$ (note that we may have $x=b$). Note that we can consider $u$ and $v$ to be the out-neighbors of $x$ on $P^b_u$ and $P^b_v$, since if those out-neighbors were resolved by any vertex $b' \in \mathcal{B}$ then $u$ and $v$ would be resolved by $b'$ too, a contradiction. Now, $u$ and $v$ cannot be in-twins, since, otherwise, one of them would be in $\mathcal{B}$, a contradiction. Thus, one vertex, say without loss of generality $u$, has an in-neighbor $w$ that is not an in-neighbor of $v$. There are now several cases to consider.
	
	If there is no arc $\overrightarrow{uw}$, then, $w$ has to be reachable from a vertex $b' \in \mathcal{B}$. However, since the underlying graph of $T$ is a tree, the only path from $b'$ to $v$ (which necessarily exists, since otherwise $b'$ resolves $u$ and $v$ since $u$ is reachable from $b'$, a contradiction) goes through $u$, and thus $\dist(b',v)=\dist(b',u)+2$, and thus $b'$ resolves $u$ and $v$, a contradiction.
	
	Hence, $w$ and $u$ are in a common strongly connected component $C$. First, assume that the arc $\overrightarrow{ux}$ does not exist, that is, $v$ is not reachable from $u$.
	Note that if any vertex from $C$ is in $\mathcal{B}$, then, $u$ and $v$ will be resolved, a contradiction.
	Now, there are only a limited number of cases where no vertex from $C$ has been added to $\mathcal{B}$. In all those cases, the underlying graph of $C$ is a path, and they are the cases that were considered neither in lines 8-17 of \Cref{alg-MDTree} nor in the special legs case.
	The first case is if $C$ is an escalator. The only possibility for this, by definition, is that $u$ and $v$ are almost-in-twins (since, otherwise, there would be another in-arc than the one at one endpoint), in which case, either $u$ or $v$ are added in $\mathcal{B}$, and thus they are resolved, a contradiction.
	The other case is if either both endpoints of $C$ have in-arcs coming from outside of $C$, or if the in-arcs coming from outside of $C$ the closest to the endpoints are not followed by any out-arc from the same vertex. However, note that, in this case, there is at least one vertex $b' \in \mathcal{B}$ such that $u$ is reachable from $b'$ but $v$ is not (there is necessarily at least one vertex in $\mathcal{B}$ "behind" every in-arc of a strongly connected component), a contradiction.
	
	This implies that the arc $\overrightarrow{ux}$ exists, and thus $x$ also belongs to $C$. First, if any vertex $y \in C$ such that $\dist(y,u) \neq \dist(y,v)$ is in $\mathcal{B}$, then, $u$ and $v$ are resolved, a contradiction. Hence, there are only two possibilities: either no vertex from $C$ is in $\mathcal{B}$, or vertices in $C \cap \mathcal{B}$ are either $x$ or in a part of $C$ that can only reach $u$ and $v$ through $x$.
	As in the previous case, if no vertex from $C$ is in $\mathcal{B}$, then, we reach a contradiction. Indeed, $C$ cannot be an escalator: either $x$ is an endpoint of $C$ and then $\overrightarrow{xv}$ is an out-arc that prevents $C$ from being an escalator, or it is not an endpoint and then the in-arc arriving at $x$ from $P^b_u$ ($x\neq b$ since no vertex from $C$ is in $\mathcal{B}$) prevents $C$ from being an escalator. In the other cases, the underlying graph of $C$ is a path with specific properties (either both endpoints have an in-arc coming from outside, or the two in-arcs coming from the outside the closest to the endpoints do not have out-arcs leaving $C$ from the same vertex), and the in-arcs coming from outside of $C$ will allow $u$ and $v$ to be resolved, a contradiction.
	
	Thus, there is at least one vertex $b' \in C \cap \mathcal{B}$, and it can only reach $u$ and $v$ through $x$. This implies that $u$ and $w$ are in a leg of $C$ with no in-arc from outside of $C$ (since, otherwise, we would have added a vertex in $\mathcal{B}$ on this side of $C$, which would resolve $u$ and $v$, a contradiction). There are three cases. First, if there is an in-arc from outside of $C$ to $x$ and $v \not\in C$, then $u$ and $w$ are in a special leg, and thus its endpoint would be in $\mathcal{B}$ and would resolve $u$ and $v$, a contradiction. Now, if there is no in-arc from outside of $C$ to $x$ and $v \not\in C$, then either $u$ and $w$ are in a special leg (a contradiction, like before), or the underlying graph of $C$ is a path, and our construction would have added the endpoint of $C$ on the side of $u$ to $\mathcal{B}$, which would resolve $u$ and $v$, a contradiction. Finally, if $v \in C$, then, $v$ has to be in a regular leg (since, otherwise, we would have added a vertex in $\mathcal{B}$ on this side of $C$, which would resolve $u$ and $v$, a contradiction), but now either the underlying graph of $C$ is a path and it would have an in-arc or we would have added one of its endpoints to $\mathcal{B}$, resolving $u$ and $v$, a contradiction; or $u$ and $v$ are the first vertices in two regular legs of $C$ spanning from the same vertex $x$ and we would have added one of the endpoints of the legs to $\mathcal{B}$, resolving $u$ and $v$, a contradiction.
	
	\ \\
	
	Thus, by our construction, every pair of vertices is resolved, and thus $\mathcal{B}$ is a resolving set. This proves that $\mathcal{B}$ is a metric basis, and thus that \Cref{alg-MDTree} is correct.
	
	Finally, it is easy to see that \Cref{alg-MDTree} computes $\mathcal{B}$ in linear time, which proves the statement of \Cref{thm-MDTree}.
\end{proof}

In~\cite{ABCHONS20}, the authors used the notion of \emph{removable source} to characterize orientations of trees with a weak metric dimension\footnote{In which one vertex may be unreachable from any resolving set vertex.} lower than the metric dimension. In the case of di-trees, however, the definition of removable source is not so clear-cut. Indeed, there are several cases where a source meets (resp. does not meet) the conditions of the removable source as defined in~\cite{ABCHONS20} and yet cannot (resp. can) be removed from a metric basis in order to obtain a weak metric basis. While we do not have a proper definition of a removable source in the context of di-trees, we get the following result:

\begin{proposition}
	\label{prop-weakMDTrees}
	There is a polynomial-time algorithm computing a weak metric basis of a di-tree.
\end{proposition}

\begin{proof}
	The result comes from two facts. The first is that every vertex we added to the metric basis in~\Cref{alg-MDTree} was necessary to either guarantee that every vertex is reached from the basis (sources) or "locally" resolve some pairs of vertices (almost-in-twins, special legs...), which still need to be resolved, and thus we cannot avoid adding those second ones to the weak metric basis either. The second is that the only possible \emph{infinite-vertex} (that is, a vertex that is not reachable from any vertex in the basis) in a directed graph is a source (\textbf{Proposition~2.2} in~\cite{ABCHONS20}). Indeed, if the infinite-vertex $s$ is not a source, then, any vertex $u$ such that there is a path from $u$ to $s$ cannot be reached from any vertex in the basis, and thus we would have several infinite-vertices, a contradiction.
	
	Hence, in a di-tree, the only possible way to have a weak metric basis is to remove a source from a metric basis without creating a pair of non-resolved vertices. This is possible in polynomial time, since the actualization of distance vectors in a di-tree will take linear time.
\end{proof}

\section{Orientations of unicyclic graphs}
\label{section-MDUnicyclic}

A unicyclic graph $U$ is constituted of a cycle $C$ with vertices $c_1,\ldots,c_n$, and each vertex $c_i$ is the root of a tree $T_i$ (we can have $T_i$ be simply the isolated $c_i$ itself).
The metric dimension of an undirected unicyclic graph has been studied in~\cite{PoissonUnicyclic,SedlarBounds,SedlarUnicyclic}. In~\cite{PoissonUnicyclic}, Poisson and Zhang proved bounds for the metric dimension of a unicyclic graph in terms of the metric dimension of a tree obtained by removing one edge from the cycle. Sedlar and \v{S}krekovski showed more recently that the metric dimension of a unicyclic graph is one of two values in~\cite{SedlarBounds}, and then the exact value of the metric dimension based on the structure of the graph in~\cite{SedlarUnicyclic}.
In this section, we will show that one can compute a metric basis of an orientation of a unicyclic graph in linear time. The algorithm mostly consists in using sources and in-twins, with a few specific edge cases to consider.

In this section, an \emph{induced directed path} $\overrightarrow{P}$ is the orientation of an induced path with only one source and one sink which are its two endpoints. It is said to be \emph{spanning from $u$} if $u$ is its source endpoint, and its \emph{length} is its number of edges. We also need the following definition:

\begin{definition}
	\label{def-badPaths}
	Let $\overrightarrow{U}$ be the orientation of a unicyclic graph. Given an orientation of a cycle $\overrightarrow{C}$ of even length $n=2k$ with two sources, if its sources are $c_i$ and $c_{i+2}$, its sinks are $c_{i+1}$ and $c_{i+1+k}$, and there are, in $\overrightarrow{C} \setminus \{c_i,c_{i+2}\}$, neither in-twins nor in-arcs coming from outside of $\overrightarrow{C}$, we call an induced directed path $\overrightarrow{P}$ a \emph{\annPath} if it verifies the three following properties:
	\begin{enumerate}
		\item $\overrightarrow{P}$ spans from $c_{i+1}$;
		\item $\overrightarrow{P}$ has length $k-2$;
		\item $\overrightarrow{P}$ has no in-arc coming from outside of $\overrightarrow{P} \cup \overrightarrow{C}$.
	\end{enumerate}
	Furthermore, if, for every vertex in $\overrightarrow{P}$ belonging to a non-empty set $I$ of in-twins, every vertex in $I$ belongs to a \annPath, then, we call $\overrightarrow{P}$ an \emph{\badPath}.
	
	A path that is a \annPath, but not an \badPath, will be called a \emph{\notSoBadPath}.
	
	Finally, a vertex might belong both to an \badPath and to a \notSoBadPath; in this case, the \notSoBadPath takes precedence (\emph{i.e.}, we will consider that the vertex belongs to the \notSoBadPath).
\end{definition}

\begin{algorithm}[!h]
	\caption{An algorithm computing the metric basis of an orientation of a unicyclic graph.}\label{alg-MDUnicyclic}
	\SetKwInOut{KwIn}{Input}
	\SetKwInOut{KwOut}{Output}
	\KwIn{An orientation $\overrightarrow{U}$ of a unicyclic graph $U$.}
	\KwOut{A metric basis $\mathcal{B}$ of $\overrightarrow{U}$.}
	
	Add to $\mathcal{B}$ every source of $\overrightarrow{U}$
	
	Apply the \textbf{special cases} in \Cref{alg-MDUnicyclic-specialCases}
	
	\ForEach{set $I$ of in-twins in $\overrightarrow{U}$ that are not already in $\mathcal{B}$}{
		\eIf{all the vertices of $I$ are in \annPaths}{
			Add $|I|-1$ vertices of $I$ to $\mathcal{B}$, prioritizing vertices in \badPaths
		}
		{
			Add $|I|-1$ vertices of $I$ to $\mathcal{B}$, prioritizing vertices in the cycle $\overrightarrow{C}$ or in \annPaths, if there are any
		}
	}
	
	\Return $\mathcal{B}$
\end{algorithm}

\begin{figure}
	\centering
	
	\begin{subfigure}[b]{0.325\linewidth}
		\centering
		\scalebox{0.75}{
			\begin{tikzpicture}
				\foreach \I in {0,...,4}{
					\node[noeud] (\I) at (360*\I/5:1) {};
				}
				\node[noeud] (a) at (2,0) {};
				\draw[\nicearrow] (a)to(0);
				\draw[\nicearrow] (0)to(1);
				\draw[\nicearrow] (1)to(2);
				\draw[\nicearrow] (2)to(3);
				\draw[\nicearrow] (3)to(4);
				\draw[\nicearrow] (4)to(0);
				
				\draw[dashed] (a) circle (0.25);
				\node at (a) [above=1.5mm] {source};
				
				\node[noeud,fill=nicered] (a) at (2,0) {};
			\end{tikzpicture}
		}
		\caption{No source in the cycle, but there is an in-arc and we are not in the second special case.}
		\label{fig-MDUnicyclicStandardCases-noSourceButInArc}
	\end{subfigure}
	\hfill
	\begin{subfigure}[b]{0.325\linewidth}
		\centering
		\scalebox{0.75}{
			\begin{tikzpicture}
				\foreach \I in {0,...,4}{
					\node[noeud] (\I) at (360*\I/5:1) {};
				}
				\begin{scope}[shift={(2,0)}]
					\node[noeud] (a) at (108:1) {};
				\end{scope}
				\draw[\nicearrow] (0)to(a);
				\draw[\nicearrow] (0)to(1);
				\draw[\nicearrow] (1)to(2);
				\draw[\nicearrow] (2)to(3);
				\draw[\nicearrow] (3)to(4);
				\draw[\nicearrow] (4)to(0);
				
				\draw[dashed,rounded corners] (0,0.6875) rectangle (2,1.1875);
				\draw (1,1.5) node {in-twins};
				\node [diamond,draw,very thick,fit=(1),inner sep=1.25] {};
				
				\node[noeud,fill=nicered] (1) at (360/5:1) {};
			\end{tikzpicture}
		}
		\caption{No source in the cycle, but there is an in-twin, which has priority.}
		\label{fig-MDUnicyclicStandardCases-noSourceButInTwin}
	\end{subfigure}
	\hfill
	\begin{subfigure}[b]{0.325\linewidth}
		\centering
		\scalebox{0.75}{
			\begin{tikzpicture}
				\foreach \I in {0,...,4}{
					\node[noeud] (\I) at (360*\I/5:1) {};
				}
				\node[noeud] (a) at (2,0) {};
				\draw[\nicearrow] (a)to(0);
				\draw[\nicearrow] (0)to(1);
				\draw[\nicearrow] (1)to(2);
				\draw[\nicearrow] (3)to(2);
				\draw[\nicearrow] (3)to(4);
				\draw[\nicearrow] (4)to(0);
				
				\draw[dashed] (a) circle (0.25);
				\node at (a) [above=1.5mm] {source};
				\draw[dashed] (3) circle (0.25);
				\node at (3) [below=1.5mm] {source};
				
				\node[noeud,fill=nicered] (a) at (2,0) {};
				\node[noeud,fill=nicered] (3) at (360*3/5:1) {};
			\end{tikzpicture}
		}
		\caption{A source $c_i$, a sink either $c_{i-1}$ or $c_{i+1}$, but an in-arc coming from outside.}
		\label{fig-MDUnicyclicStandardCases-OneSinkAndInArc}
	\end{subfigure}
	\vskip\baselineskip
	\begin{subfigure}[b]{0.325\linewidth}
		\centering
		\scalebox{0.75}{
			\begin{tikzpicture}
				\foreach \I in {0,...,4}{
					\node[noeud] (\I) at (360*\I/5:1) {};
				}
				\draw[\nicearrow] (1)to(0);
				\draw[\nicearrow] (1)to(2);
				\draw[\nicearrow] (2)to(3);
				\draw[\nicearrow] (3)to(4);
				\draw[\nicearrow] (0)to(4);
				
				\draw[dashed] (1) circle (0.25);
				\node at (1) [above=1.5mm] {source};
				\coordinate (2B) at ($(2) + (-0.375,-0.1875)$);
				\coordinate (0B) at ($(0) + (0.375,0.1875)$);
				\draw[rounded corners,dashed,rotate around={-17.5:(144:1)}] (2B) rectangle (0B);
				\draw (0.125,0.3125) node {\rotatebox{-17.5}{in-twins}};
				
				\node[noeud,fill=nicered] (1) at (360/5:1) {};
				\node[noeud,fill=nicered] (2) at (360*2/5:1) {};
			\end{tikzpicture}
		}
		\caption{A source $c_i$, and a sink that is neither $c_{i-1}$ nor $c_{i+1}$.}
		\label{fig-MDUnicyclicStandardCases-OneSinkFarFromTheSource}
	\end{subfigure}
	\hfill
	\begin{subfigure}[b]{0.325\linewidth}
		\centering
		\scalebox{0.75}{
			\begin{tikzpicture}
				\foreach \I in {0,...,5}{
					\node[noeud] (\I) at (360*\I/6:1) {};
				}
				\begin{scope}[shift={(2,0)}]
					\node[noeud] (a) at (120:1) {};
					\node[noeud] (b) at (240:1) {};
				\end{scope}
				\node[noeud] (c) at (-0.5,0) {};
				\draw[\nicearrow] (1)to(0);
				\draw[\nicearrow] (2)to(1);
				\draw[\nicearrow] (3)to(2);
				\draw[\nicearrow] (3)to(4);
				\draw[\nicearrow] (4)to(5);
				\draw[\nicearrow] (5)to(0);
				\draw[\nicearrow] (1)to(a);
				\draw[\nicearrow] (5)to(b);
				\draw[\nicearrow] (3)to(c);
				
				\draw[dashed] (3) circle (0.25);
				\node at (3) [left=2mm] {source};
				\draw[dashed,rounded corners] (237.5:1.25) rectangle (107.5:1.125);
				\node at (2) [above=1.5mm] {in-twins};
				\node at (0) [right=1.5mm] {$c_{i+k}$};
				\node at (a) [right=1mm] {$v_-$};
				\node at (b) [right=1mm] {$v_+$};
				
				\node[noeud,fill=nicered] (2) at (120:1) {};
				\node[noeud,fill=nicered] (3) at (180:1) {};
				\node[noeud,fill=nicered] (4) at (240:1) {};
			\end{tikzpicture}
		}
		\caption{$|\overrightarrow{C}|=2k$, there is one source $c_i$ in $\overrightarrow{C}$, the sink is $c_{i+k}$, $c_{i+k-1}$ and $c_{i+k+1}$ both have out-neighbors $v_-$ and $v_+$ with in-degree~1, no vertex in $\overrightarrow{C} \setminus \{c_{i}\}$ has an in-arc, no vertex in $\overrightarrow{C} \setminus \{c_{i},c_{i-1},c_{i+1}\}$ is an in-twin, but $c_{i-1}$ and $c_{i+1}$ are in a set of in-twins of size at least~3.}
		\label{fig-MDUnicyclicStandardCases-oneSinkCIPlus2K}
	\end{subfigure}
	\hfill
	\begin{subfigure}[b]{0.325\linewidth}
		\centering
		\scalebox{0.75}{
			\begin{tikzpicture}
				\foreach \I in {0,...,5}{
					\node[noeud] (\I) at (360*\I/6:1) {};
				}
				\node[noeud] (a) at (2,0.5) {};
				\node[noeud] (b) at (2,-0.5) {};
				\coordinate (2B) at ($(2) + (-1,0)$);
				\node[noeud] (c) at (2B) {};
				\draw[\nicearrow,line width=0.5mm] (0)to(a);
				\draw[\nicearrow,line width=0.5mm] (0)to(b);
				\draw[\nicearrow] (c)to(2);
				\draw[\nicearrow] (1)to(0);
				\draw[\nicearrow] (1)to(2);
				\draw[\nicearrow] (2)to(3);
				\draw[\nicearrow] (4)to(3);
				\draw[\nicearrow] (5)to(4);
				\draw[\nicearrow] (5)to(0);
				
				\draw[dashed] (1) circle (0.25);
				\node at (1) [above=1.5mm] {source};
				\draw[dashed] (5) circle (0.25);
				\node at (5) [below=1.5mm] {source};
				\draw[dashed] (c) circle (0.25);
				\node at (c) [above=1.5mm] {source};
				\draw[dashed,rounded corners] (1.75,-0.75) rectangle (2.25,0.75);
				\draw (2,1) node {in-twins};
				
				\node[noeud,fill=nicered] (a) at (2,0.5) {};
				\node[noeud,fill=nicered] (c) at (2B) {};
				\node[noeud,fill=nicered] (1) at (360/6:1) {};
				\node[noeud,fill=nicered] (5) at (360*5/6:1) {};
			\end{tikzpicture}
		}
		\caption{There are two sources $c_i$ and $c_{i+2}$ in $\overrightarrow{C}$, the sinks are $c_{i+1}$ and $c_{i+1+\frac{n}{2}}$, and there is at least one \badPath (in bolded arcs) and no \notSoBadPath, but there is an in-arc coming from outside of $\overrightarrow{C}$.}
		\label{fig-MDUnicyclicStandardCases-allBadPathsButInArc}
	\end{subfigure}
	\vskip\baselineskip
	\begin{subfigure}[b]{0.325\linewidth}
		\centering
		\scalebox{0.75}{
			\begin{tikzpicture}
				\foreach \I in {0,...,9}{
					\node[noeud] (\I) at (36*\I:1) {};
				}
				\node[noeud] (a0) at (2,0.75) {};
				\node[noeud] (a1) at (3,0.75) {};
				\node[noeud] (a2) at (4,0.75) {};
				\node[noeud] (b0) at (2,0) {};
				\node[noeud] (b1) at (3,0) {};
				\node[noeud] (b2) at (4,0) {};
				\node[noeud] (c) at (3,-0.75) {};
				\draw[\nicearrow,line width=0.5mm] (0)to(a0);
				\draw[\nicearrow,line width=0.5mm] (a0)to(a1);
				\draw[\nicearrow,line width=0.5mm] (a1)to(a2);
				\draw[\nicearrow,dashed,line width=0.5mm] (0)to(b0);
				\draw[\nicearrow,dashed,line width=0.5mm] (b0)to(b1);
				\draw[\nicearrow,dashed,line width=0.5mm] (b1)to(b2);
				\draw[\nicearrow] (b0)to(c);
				\draw[\nicearrow] (1)to(0);
				\draw[\nicearrow] (1)to(2);
				\draw[\nicearrow] (2)to(3);
				\draw[\nicearrow] (3)to(4);
				\draw[\nicearrow] (4)to(5);
				\draw[\nicearrow] (6)to(5);
				\draw[\nicearrow] (7)to(6);
				\draw[\nicearrow] (8)to(7);
				\draw[\nicearrow] (9)to(8);
				\draw[\nicearrow] (9)to(0);
				
				\draw[dashed] (1) circle (0.25);
				\node [below left=1mm and -1mm of 1] {source};
				\draw[dashed] (9) circle (0.25);
				\node [above left=1mm and -1mm of 9] {source};
				\draw[dashed,rounded corners] (1.75,-0.25) rectangle (2.25,1);
				\draw (2,1.25) node {in-twins};
				\draw[dashed,rounded corners] (2.75,0.25) rectangle (3.25,-1);
				\draw (3,-1.25) node {in-twins};
				\node [diamond,draw,very thick,fit=(a0),inner sep=1.25] {};
				\node [diamond,draw,very thick,fit=(b1),inner sep=1.25] {};
				
				\node[noeud,fill=nicered] (1) at (36:1) {};
				\node[noeud,fill=nicered] (9) at (36*9:1) {};
				\node[noeud,fill=nicered] (a0) at (2,0.75) {};
				\node[noeud,fill=nicered] (b1) at (3,0) {};
			\end{tikzpicture}
		}
		\caption{There are two sources $c_i$ and $c_{i+2}$ in $\overrightarrow{C}$, the sinks are $c_{i+1}$ and $c_{i+1+\frac{n}{2}}$, there is neither in-twin nor in-arc in $\overrightarrow{C} \setminus \{c_i,c_{i+2}\}$, and there is at least one \badPath (in bolded arcs), but there is at least one \notSoBadPath (in bolded and dashed arcs).}
		\label{fig-MDUnicyclicStandardCases-badPathsAndNotSoBadPaths}
	\end{subfigure}
	\hfill
	\begin{subfigure}[b]{0.325\linewidth}
		\centering
		\scalebox{0.75}{
			\begin{tikzpicture}
				\foreach \I in {0,...,5}{
					\node[noeud] (\I) at (360*\I/6:1) {};
				}
				\node[noeud] (a) at (2,0) {};
				\node[noeud] (b) at (-2,0) {};
				\draw[\nicearrow] (0)to(a);
				\draw[\nicearrow] (3)to(b);
				\draw[\nicearrow] (1)to(0);
				\draw[\nicearrow] (1)to(2);
				\draw[\nicearrow] (2)to(3);
				\draw[\nicearrow] (3)to(4);
				\draw[\nicearrow] (5)to(4);
				\draw[\nicearrow] (5)to(0);
				
				\draw[dashed] (1) circle (0.25);
				\node [above=1mm of 1] {source};
				\draw[dashed] (5) circle (0.25);
				\node [below=1mm of 5] {source};
				
				\node[noeud,fill=nicered] (1) at (60:1) {};
				\node[noeud,fill=nicered] (5) at (300:1) {};
			\end{tikzpicture}
		}
		\caption{Two sources, not in the special case configuration.}
		\label{fig-MDUnicyclicStandardCases-twoSources}
	\end{subfigure}
	\hfill
	\begin{subfigure}[b]{0.325\linewidth}
		\centering
		\scalebox{0.75}{
			\begin{tikzpicture}
				\foreach \I in {0,...,7}{
					\node[noeud] (\I) at (360*\I/8:1) {};
				}
				\node[noeud] (a) at (2,0) {};
				\node[noeud] (b) at (-2,0) {};
				\draw[\nicearrow] (0)to(a);
				\draw[\nicearrow] (4)to(b);
				\draw[\nicearrow] (1)to(0);
				\draw[\nicearrow] (2)to(1);
				\draw[\nicearrow] (2)to(3);
				\draw[\nicearrow] (3)to(4);
				\draw[\nicearrow] (5)to(4);
				\draw[\nicearrow] (5)to(6);
				\draw[\nicearrow] (7)to(6);
				\draw[\nicearrow] (7)to(0);
				
				\draw[dashed] (2) circle (0.25);
				\node [above=1mm of 2] {source};
				\draw[dashed] (5) circle (0.25);
				\node [below=1mm of 5] {source};
				\draw[dashed] (7) circle (0.25);
				\node [below=1mm of 7] {source};
				
				\node[noeud,fill=nicered] (2) at (90:1) {};
				\node[noeud,fill=nicered] (5) at (225:1) {};
				\node[noeud,fill=nicered] (7) at (315:1) {};
			\end{tikzpicture}
		}
		\caption{More than two sources.}
		\label{fig-MDUnicyclicStandardCases-moreThanTwoSources}
	\end{subfigure}
	
	\caption{The standard cases of \Cref{alg-MDUnicyclic}, when a metric basis $\mathcal{B}$ of an orientation of a unicyclic graph contains every source and resolves every set of in-twins (with some priority, identified with a diamond). Vertices in $\mathcal{B}$ are in red.}
	\label{fig-MDUnicyclicStandardCases}
\end{figure}
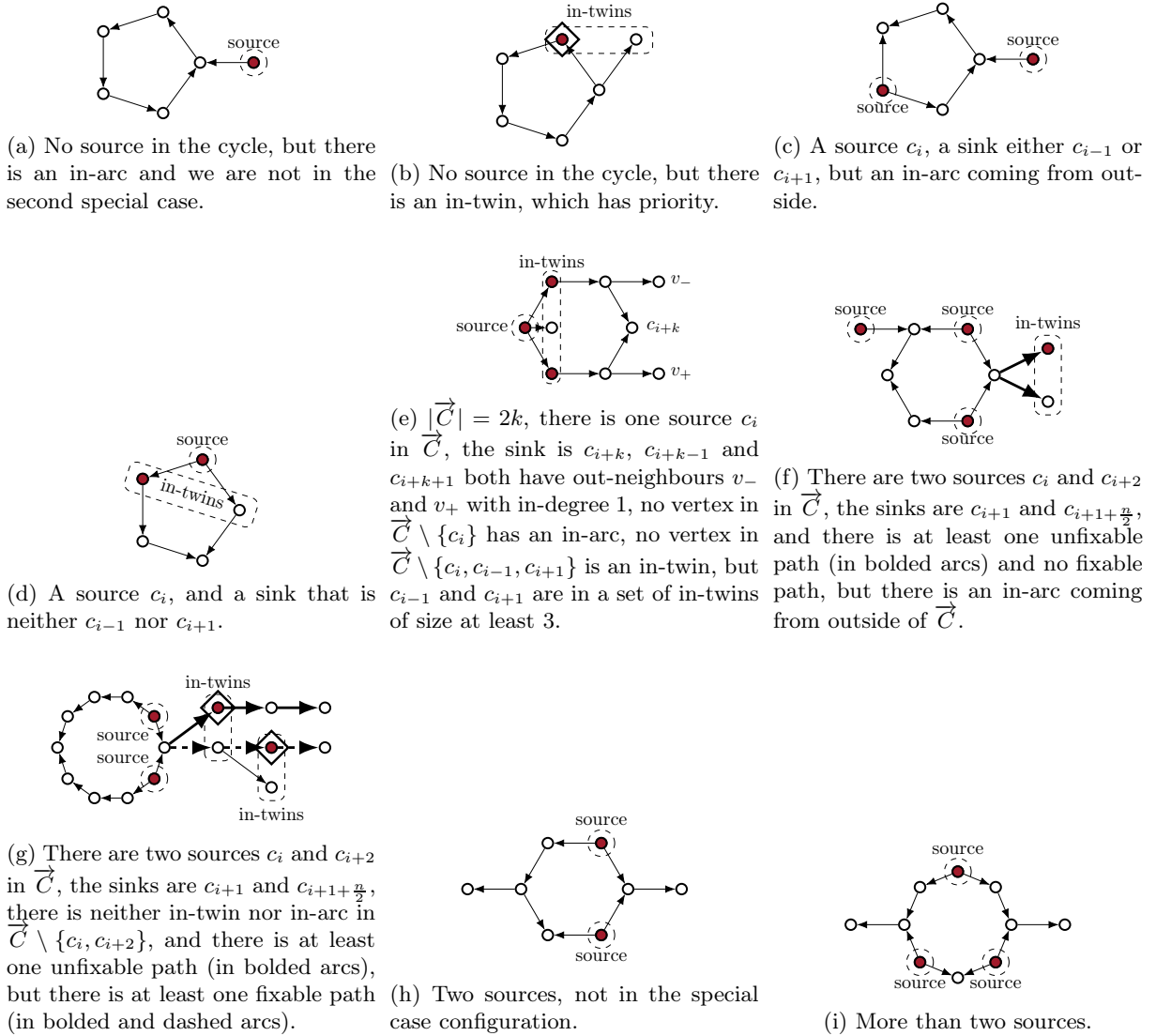

\begin{algorithm}[!h]
	\caption{Special cases of \Cref{alg-MDUnicyclic}.}\label{alg-MDUnicyclic-specialCases}
	
	\If{the cycle $\overrightarrow{C}$ has no sink, there is no in-arc coming from outside of $C$, and no vertex of $C$ is in a set of in-twin}{
		Add $c_1$ to $\mathcal{B}$
	}
	
	\If{the cycle $\overrightarrow{C}$ has no sink, there is exactly one in-arc $\overrightarrow{uc_i}$ with $u \not\in \overrightarrow{C}$, no vertex $c_j$ with $j \neq i$ is an in-twin or has in-arc coming from outside of $\overrightarrow{C}$, and $u$ has an out-neighbor $v$ with $N^-(v)=\{u\}$}{
		Add $c_i$ to $\mathcal{B}$
	}
	
	\If{the cycle $\overrightarrow{C}$ has exactly one source $c_i$}{
		\If{the one sink is either $c_{i-1}$ or $c_{i+1}$, and no vertex $c_j$ with $j \neq i$ is an in-twin or has an in-arc coming from outside of $\overrightarrow{C}$}{
			Add $c_{i-1}$ to $\mathcal{B}$
		}
		\ElseIf{the one sink is $c_{i+k}$ with $k>1$, $|\overrightarrow{C}|\geq 2k$, $c_{i+k-1}$ (resp. $c_{i+k+1}$) has an out-neighbor $v$ such that $N^-(v)=\{c_{i+k-1}\}$ (resp. $N^-(v)=\{c_{i+k+1}\}$), no vertex in $\{c_{i-1},c_{i-2},\ldots,c_{i+k}\}$ (resp. $\{c_{i+1},c_{i+2},\ldots,c_{i+k}\}$) has an in-arc, and no vertex in $\{c_{i-2},c_{i-3},\ldots,c_{i+k+1}\}$ (resp. $\{c_{i+2},c_{i+3},\ldots,c_{i+k-1}\}$) is an in-twin}{
			Add $c_{i-1}$ (resp. $c_{i+1}$) to $\mathcal{B}$
		}
		\ElseIf{the one sink is $c_{i+k}$ with $k>1$, $|\overrightarrow{C}|=2k$, $c_{i+k-1}$ has an out-neighbor $v_-$ such that $N^-(v_-)=\{c_{i+k-1}\}$, $c_{i+k+1}$ has an out-neighbor $v_+$ such that $N^-(v_+)=\{c_{i+k+1}\}$, no vertex in $\overrightarrow{C}$ except $c_i$ has an in-arc, no vertex in $\overrightarrow{C} \setminus \{c_i, c_{i-1}, c_{i+1} \}$ is an in-twin, and $c_{i-1}$ and $c_{i+1}$ are not in a set $I$ of in-twins verifying $|I|\geq 3$}{
			Add $c_{i+k}$ to $\mathcal{B}$
		}
	}
	
	\If{the cycle $\overrightarrow{C}$ has exactly two sources $c_i$ and $c_{i+2}$, $|\overrightarrow{C}|=2k$ with $k>2$, the two sinks are $c_{i+1}$ and $c_{i+1+k}$, no vertex from $\overrightarrow{C}$ except $c_i$ and $c_{i+2}$ is an in-twin or has an in-arc coming from outside of $\overrightarrow{C}$, there is at least one \badPath, and there is no \notSoBadPath}{
		Add $c_{i+1}$ to $\mathcal{B}$
	}
	
\end{algorithm}

\ \\\noindent\textbf{Explanation of \Cref{alg-MDUnicyclic}.} The algorithm will compute a metric basis $\mathcal{B}$ of an orientation $\overrightarrow{U}$ of a unicyclic graph $U$ in linear time. The result on several cases is depicted in \Cref{fig-MDUnicyclicSpecialCases,fig-MDUnicyclicStandardCases}.
The first thing we do is to add every source in $\overrightarrow{U}$ to $\mathcal{B}$ (line~1). We will also manage the sets of in-twins in $\overrightarrow{U}$ (lines 3-7), which we need to do after taking care of some special cases that might influence the choice of in-twins. When we have the choice, we prioritize taking in-twins that are in the cycle to guarantee reachability of vertices in the cycle. Note that those two sets (along with the right priority) are enough in most cases, as depicted in \Cref{fig-MDUnicyclicStandardCases}.

We then have to manage six specific cases (line 2). \textbf{Those special cases are handled in \Cref{alg-MDUnicyclic-specialCases}, to which the line numbers in the next five paragraphs will refer.} The first two special cases occur when the cycle has no sink. First, if the cycle has no sink, no in-twin, and no arc coming from outside, then, we have to add one vertex of the cycle to $\mathcal{B}$ in order to maintain reachability (lines 1-2, depicted in \Cref{fig-MDUnicyclicSpecialCases-nothingInC}). Then, if the cycle has no sink, only one in-arc $\overrightarrow{uc_i}$ is coming from outside of it, and there is a vertex $v$ with $N^-(v)=\{u\}$, then, we add $c_i$ or $v$ to $\mathcal{B}$ in order to resolve them (lines 3-4, depicted in \Cref{fig-MDUnicyclicSpecialCases-noSourceAndSpecificInArc}).

The next three special cases occur when the cycle has one sink. First, if there is only one sink in the cycle, it is an out-neighbor of the source, and no vertex from the cycle apart from the source is an in-twin or has an in-arc coming from outside of the cycle, then we need to add one of the out-neighbors of the source in the cycle to $\mathcal{B}$ in order to resolve them (lines 6-7), depicted in \Cref{fig-MDUnicyclicSpecialCases-oneSinkAndNothingElseInC}).

Then, there are two specific cases when the cycle has one sink, both based on the same principle. Both happen when the source is $c_i$, the sink is $c_{i+k}$, it has no in-arc, and the cycle contains at least $2k$ vertices. In the fourth special case (lines~8-9), depicted in \Cref{fig-MDUnicyclicSpecialCases-oneSinkCIPlusKLarge}),
the vertex $c_{i+k-1}$ has an out-neighbor $v$ verifying $N^-(v)=\{c_{i+k-1}\}$. We can see that, if no vertex in the other path from $c_i$ to $c_{i+k}$ (the path going through $c_{i-1},c_{i-2},\ldots,c_{i+k+1}$) is in $\mathcal{B}$, then, $v$ and $c_{i+k}$ will not be resolved. Those vertices can be added to $\mathcal{B}$ if they have an in-arc or if they are an in-twin (they will have priority). However, note that $c_{i-1}$ might be an in-twin of $c_{i+1}$, in which case it should be added to $\mathcal{B}$, resolving the conflict. Hence, if none of $c_{i-1},c_{i-2},\ldots,c_{i+k+1}$ has an in-arc or is an in-twin, then, we can add $c_{i-1}$ to $\mathcal{B}$ in order to resolve $v$ and $c_{i+k}$. Note that, in this case, in comparison to just the sources and the resolution of sets of in-twins, we add one more vertex to $\mathcal{B}$ if $c_{i-1}$ is the only in-twin of $c_{i+1}$. The same reasoning can be made with the symmetric case.

The fifth special case (lines~10-11), depicted in \Cref{fig-MDUnicyclicSpecialCases-oneSinkCIPlus2K}),
occurs when the cycle contains exactly $2k$ vertices and both $c_{i+k-1}$ and $c_{i+k+1}$ have an out-neighbor (respectively $v_-$ and $v_+$) with in-degree~1: the pairs of vertices $(v_-,c_{i+k})$ and $(v_+,c_{i+k})$ might not be resolved. We can see that any in-arc or in-twin along a path from $c_i$ to $c_{i+k}$ will resolve $c_{i+k}$ and the pendent $v$ on the other path (thus either fully resolving those two pairs, or bringing us back to the previous special case), except if $c_{i-1}$ and $c_{i+1}$ are the only in-twins in the cycle and if they do not have another in-twin. Hence, if no vertex from the cycle except $c_i$ has an in-arc, no vertex from the cycle except $c_i$, $c_{i-1}$ and $c_{i+1}$ is an in-twin, and $c_{i-1}$ and $c_{i+1}$ do not have another in-twin, then, we need to add at least one more vertex to $\mathcal{B}$ in order to resolve the two pairs of vertices, and adding $c_{i+k}$ does exactly that.

Finally, the sixth special case is more complex (lines 12-13, depicted in \Cref{fig-MDUnicyclicSpecialCases-allBadPaths},
and consideration in the choice of in-twins, depicted in \Cref{fig-MDUnicyclicStandardCases-badPathsAndNotSoBadPaths}).
Assume that the cycle $\overrightarrow{C}$ is of even length $n$, has neither in-twin nor in-arc coming from outside (except the sources), and that there are two sinks in $\overrightarrow{C}$: one at distance~1 from the sources, and the other at the opposite end of $\overrightarrow{C}$.
Now, if the first sink has spanning \annPaths, then, the second sink and the endpoints of those \annPaths might not be resolved, since they are at the same distance ($\frac{n}{2}-1$) from both sources of $\overrightarrow{C}$. Thus, we need to apply a strategy in order to resolve those vertices while trying to not add a supplementary vertex to $\mathcal{B}$. This is done by considering the two kinds of \annPaths, and having a priority in the selection of in-twins. The details will be in the proof.

All the other cases of the cycle are already resolved through the sources and in-twins steps.

\begin{figure}
	\centering
	
	\begin{subfigure}[b]{0.325\linewidth}
		\centering
		\scalebox{0.75}{
			\begin{tikzpicture}
				\foreach \I in {0,...,4}{
					\node[noeud] (\I) at (360*\I/5:1) {};
				}
				\node[noeud] (a) at (2,0) {};
				\node[noeud] (b) at (3,0) {};
				\draw[\nicearrow] (b)to(a);
				\draw[\nicearrow] (0)to(a);
				\draw[\nicearrow] (0)to(1);
				\draw[\nicearrow] (1)to(2);
				\draw[\nicearrow] (2)to(3);
				\draw[\nicearrow] (3)to(4);
				\draw[\nicearrow] (4)to(0);
				
				\draw[dashed] (b) circle (0.25);
				\node at (b) [above=1.5mm] {source};
				\node [rectangle,draw,very thick,fit=(1),inner sep=3] {};
				\node at (1) [above=1.5mm] {$c_1$};
				
				\node[noeud,fill=nicered] (b) at (3,0) {};
				\node[noeud,fill=nicered] (1) at (360/5:1) {};
			\end{tikzpicture}
		}
		\caption{First special case: there is neither source, nor in-twin or in-arc coming from outside in $\overrightarrow{C}$.}
		\label{fig-MDUnicyclicSpecialCases-nothingInC}
	\end{subfigure}
	\hfill
	\begin{subfigure}[b]{0.325\linewidth}
		\centering
		\scalebox{0.75}{
			\begin{tikzpicture}
				\foreach \I in {0,...,4}{
					\node[noeud] (\I) at (360*\I/5:1) {};
				}
				\node[noeud] (a) at (2,0) {};
				\node[noeud] (b) at (3,0) {};
				\draw[\nicearrow] (a)to(0);
				\draw[\nicearrow] (a)to(b);
				\draw[\nicearrow] (0)to(1);
				\draw[\nicearrow] (1)to(2);
				\draw[\nicearrow] (2)to(3);
				\draw[\nicearrow] (3)to(4);
				\draw[\nicearrow] (4)to(0);
				
				\draw[dashed] (a) circle (0.25);
				\node at (a) [above=1.5mm] {source};
				\node [rectangle,draw,very thick,fit=(0),inner sep=3] {};
				\node at (0) [left=1.5mm] {$c_i$};
				
				\node[noeud,fill=nicered] (a) at (2,0) {};
				\node[noeud,fill=nicered] (0) at (0:1) {};
			\end{tikzpicture}
		}
		\caption{Second special case: there is no source or in-twin in the cycle, and there is a unique in-arc $\overrightarrow{uc_i}$ and a vertex $v$ with $N^-(v)=\{u\}$.}
		\label{fig-MDUnicyclicSpecialCases-noSourceAndSpecificInArc}
	\end{subfigure}
	\hfill
	\begin{subfigure}[b]{0.325\linewidth}
		\centering
		\scalebox{0.75}{
			\begin{tikzpicture}
				\foreach \I in {0,...,4}{
					\node[noeud] (\I) at (360*\I/5:1) {};
				}
				\node[noeud] (a) at (2,0) {};
				\draw[\nicearrow] (a)to(0);
				\draw[\nicearrow] (0)to(1);
				\draw[\nicearrow] (2)to(1);
				\draw[\nicearrow] (3)to(2);
				\draw[\nicearrow] (4)to(3);
				\draw[\nicearrow] (0)to(4);
				
				\draw[dashed] (a) circle (0.25);
				\node at (a) [above=1.5mm] {source};
				\node [rectangle,draw,very thick,fit=(1),inner sep=3] {};
				\node at (1) [above=1.5mm] {$c_{i-1}$};
				
				\node[noeud,fill=nicered] (a) at (2,0) {};
				\node[noeud,fill=nicered] (1) at (360/5:1) {};
			\end{tikzpicture}
		}
		\caption{Third special case: there is one source $c_i$ in $\overrightarrow{C}$, the sink is either $c_{i-1}$ or $c_{i+1}$, and there is neither in-twin nor in-arc in $\overrightarrow{C} \setminus \{c_i\}$.}
		\label{fig-MDUnicyclicSpecialCases-oneSinkAndNothingElseInC}
	\end{subfigure}
	\vskip\baselineskip
	\begin{subfigure}[b]{0.325\linewidth}
		\centering
		\scalebox{0.75}{
			\begin{tikzpicture}
				\foreach \I in {0,...,4}{
					\node[noeud] (\I) at (360*\I/5:1) {};
				}
				\begin{scope}[shift={(2,0)}]
					\node[noeud,fill=nicered] (a) at (108:1) {};
					\node[noeud] (b) at (252:1) {};
				\end{scope}
				\draw[\nicearrow] (a)to(0);
				\draw[\nicearrow] (0)to(b);
				\draw[\nicearrow] (1)to(0);
				\draw[\nicearrow] (1)to(2);
				\draw[\nicearrow] (2)to(3);
				\draw[\nicearrow] (3)to(4);
				\draw[\nicearrow] (0)to(4);
				
				\draw[dashed] (a) circle (0.25);
				\node at (a) [above=1.5mm] {source};
				\draw[dashed] (1) circle (0.25);
				\node at (1) [above=1.5mm] {source};
				\node [rectangle,draw,very thick,fit=(2),inner sep=3] {};
				\node at (4) [below=1mm] {$c_{i+k}$};
				\node at (b) [below=1mm] {$v$};
				
				\node[noeud,fill=nicered] (1) at (360/5:1) {};
				\node[noeud,fill=nicered] (2) at (144:1) {};
			\end{tikzpicture}
		}
		\caption{Fourth special case: $|\overrightarrow{C}|\geq 2k$, there is one source $c_i$ in $\overrightarrow{C}$, the sink is $c_{i+k}$, $c_{i+k-1}$ (w.l.o.g.) has an out-neighbor $v$ with in-degree~1, no vertex in $\{c_{i-1},c_{i-2},\ldots,c_{i+k}\}$ has an in-arc, and no vertex in $\{c_{i-2},c_{i-3},\ldots,c_{i+k+1}\}$ is an in-twin.}
		\label{fig-MDUnicyclicSpecialCases-oneSinkCIPlusKLarge}
	\end{subfigure}
	\hfill
	\begin{subfigure}[b]{0.325\linewidth}
		\centering
		\scalebox{0.75}{
			\begin{tikzpicture}
				\foreach \I in {0,...,5}{
					\node[noeud] (\I) at (360*\I/6:1) {};
				}
				\begin{scope}[shift={(2,0)}]
					\node[noeud] (a) at (120:1) {};
					\node[noeud] (b) at (240:1) {};
				\end{scope}
				\draw[\nicearrow] (1)to(0);
				\draw[\nicearrow] (2)to(1);
				\draw[\nicearrow] (3)to(2);
				\draw[\nicearrow] (3)to(4);
				\draw[\nicearrow] (4)to(5);
				\draw[\nicearrow] (5)to(0);
				\draw[\nicearrow] (1)to(a);
				\draw[\nicearrow] (5)to(b);
				
				\draw[dashed] (3) circle (0.25);
				\node at (3) [left=2mm] {source};
				\draw[dashed,rounded corners] (237.5:1.25) rectangle (107.5:1.125);
				\node at (2) [above=1.5mm] {in-twins};
				\node [rectangle,draw,very thick,fit=(0),inner sep=3] {};
				\node at (0) [right=1.5mm] {$c_{i+k}$};
				\node at (a) [right=1mm] {$v_-$};
				\node at (b) [right=1mm] {$v_+$};
				
				\node[noeud,fill=nicered] (3) at (180:1) {};
				\node[noeud,fill=nicered] (4) at (240:1) {};
				\node[noeud,fill=nicered] (0) at (0:1) {};
			\end{tikzpicture}
		}
		\caption{Fifth special case: $|\overrightarrow{C}|=2k$, there is one source $c_i$ in $\overrightarrow{C}$, the sink is $c_{i+k}$, $c_{i+k-1}$ and $c_{i+k+1}$ both have out-neighbors $v_-$ and $v_+$ with in-degree~1, no vertex in $\overrightarrow{C} \setminus \{c_{i}\}$ has an in-arc, no vertex in $\overrightarrow{C} \setminus \{c_{i},c_{i-1},c_{i+1}\}$ is an in-twin, and $c_{i-1}$ and $c_{i+1}$ are not in a set of in-twins of size at least~3.}
		\label{fig-MDUnicyclicSpecialCases-oneSinkCIPlus2K}
	\end{subfigure}
	\hfill
	\begin{subfigure}[b]{0.325\linewidth}
		\centering
		\scalebox{0.75}{
			\begin{tikzpicture}
				\foreach \I in {0,...,5}{
					\node[noeud] (\I) at (360*\I/6:1) {};
				}
				\node[noeud] (a) at (2,0.5) {};
				\node[noeud] (b) at (2,-0.5) {};
				\node[noeud] (c) at (3,-0.5) {};
				\draw[\nicearrow,line width=0.5mm] (0)to(a);
				\draw[\nicearrow,line width=0.5mm] (0)to(b);
				\draw[\nicearrow] (b)to(c);
				\draw[\nicearrow] (1)to(0);
				\draw[\nicearrow] (1)to(2);
				\draw[\nicearrow] (2)to(3);
				\draw[\nicearrow] (4)to(3);
				\draw[\nicearrow] (5)to(4);
				\draw[\nicearrow] (5)to(0);
				
				\draw[dashed] (1) circle (0.25);
				\node at (1) [above=1.5mm] {source};
				\draw[dashed] (5) circle (0.25);
				\node at (5) [below=1.5mm] {source};
				\draw[dashed,rounded corners] (1.75,-0.75) rectangle (2.25,0.75);
				\draw (2,1) node {in-twins};
				\node [rectangle,draw,very thick,fit=(0),inner sep=3] {};
				\node at (0) [left=1.5mm] {$c_{i+1}$};
				
				\node[noeud,fill=nicered] (a) at (2,0.5) {};
				\node[noeud,fill=nicered] (0) at (0:1) {};
				\node[noeud,fill=nicered] (1) at (360/6:1) {};
				\node[noeud,fill=nicered] (5) at (360*5/6:1) {};
			\end{tikzpicture}
		}
		\caption{Sixth special case: there are two sources $c_i$ and $c_{i+2}$ in $\overrightarrow{C}$, the sinks are $c_{i+1}$ and $c_{i+1+\frac{n}{2}}$, there is neither in-twin nor in-arc in $\overrightarrow{C} \setminus \{c_i,c_{i+2}\}$, and there is at least one \badPath (in bolded arcs) and no \notSoBadPath.}
		\label{fig-MDUnicyclicSpecialCases-allBadPaths}
	\end{subfigure}
	
	\caption{The special cases of \Cref{alg-MDUnicyclic}, managed in \Cref{alg-MDUnicyclic-specialCases}, where we have to add to $\mathcal{B}$ one more vertex other than the sources and the resolution of sets of in-twins. Vertices in $\mathcal{B}$ are in red, and the supplementary vertex added to $\mathcal{B}$ is identified by a square around it.}
	\label{fig-MDUnicyclicSpecialCases}
\end{figure}
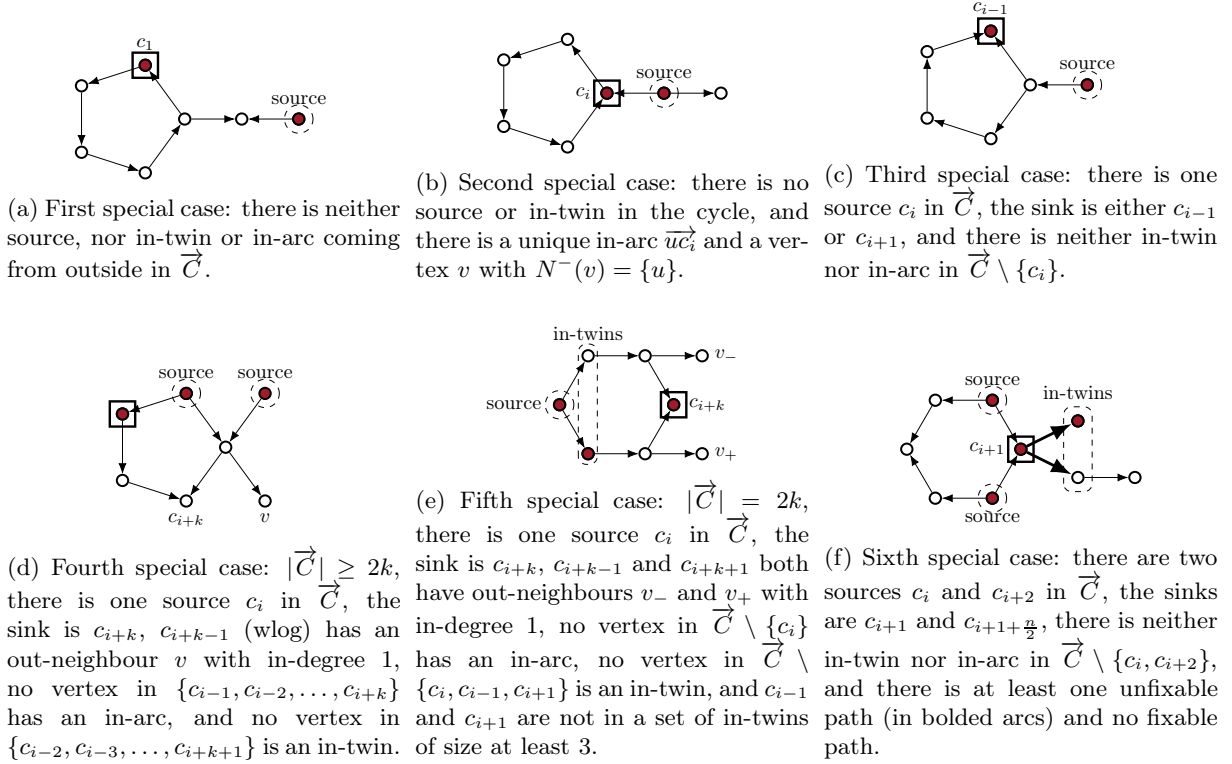

\begin{theorem}
	\label{thm-MDUnicyclic}
	\Cref{alg-MDUnicyclic} computes a metric basis of an orientation of a unicyclic graph in linear time.
\end{theorem}

\begin{proof}
	Let $\overrightarrow{U}$ be an orientation of a unicyclic graph with cycle $\overrightarrow{C}$.
	First, note that sources and $|I|-1|$ of each set $I$ of in-twins must be in the metric basis $\mathcal{B}$. Furthermore, when a vertex of a set $I$ of in-twins is in $\overrightarrow{C}$, then, it should be prioritirized for reachability reasons. However, this is not sufficient to obtain a metric basis, since either some vertices may be unreachable or some pairs of vertices may be non-resolved, which is why we will need the six special cases. In those cases, we will either have to give stronger priority to some in-twins, or have to add one more vertex to $\mathcal{B}$.
	
	The first case is if vertices in $\overrightarrow{C}$ are not reachable from any source or in-twin. This is only possible if $\overrightarrow{C}$ contains no sink, no in-arc is coming from outside of $\overrightarrow{C}$, and no vertex of $\overrightarrow{C}$ is an in-twin. In this case, we need to add any vertex from $\overrightarrow{C}$ to $\mathcal{B}$.
	
	The second case is if there is exactly one in-arc $\overrightarrow{uc_i}$ with $u \not\in \overrightarrow{C}$, the cycle has no sink, no $c_j$ with $j \neq i$ is an in-twin or has an in-arc coming from outside of $\overrightarrow{C}$, and there is a vertex $v$ verifying $N^-(v)=\{u\}$. In this case, $v$ and $c_i$ are not resolved, and thus we need to add at least one of them to $\mathcal{B}$.
	
	The third case is if there is only one source $c_i$ in $\overrightarrow{C}$, the sink is either $c_{i-1}$ or $c_{i+1}$, and no vertex $c_j$ with $j \neq i$ is an in-twin or has any in-arc coming from outside of $\overrightarrow{C}$. In this case, the out-neighbors of $c_i$ are not in-twins but cannot be resolved without taking at least one of them into $\mathcal{B}$.
	
	The fourth and fifth case are linked. In both cases, the cycle $\overrightarrow{C}$ contains at least $2k$ vertices, one source $c_i$, and one sink $c_{i+k}$ with no in-arc. The problem will be when an in-neighbor of the sink has an out-neighbor $v$ with in-degree~1: $v$ and $c_{i+k}$ might not be resolved. Let us see when this can happen.
	
	In the fourth case, either $c_{i+k-1}$ or $c_{i+k+1}$ has one out-neighbor $v$ with in-degree~1. Without loss of generality, we will consider that it is $c_{i+k-1}$. Assume furthermore that no vertex from $\{c_{i-1},c_{i-2},\ldots,c_{i+k+1}\}$ has an in-arc, and no vertex from $\{c_{i-2},c_{i-3},\ldots,c_{i+k+1}\}$ is an in-twin (since, otherwise, $v$ and $c_{i+k}$ would be resolved). Now, whether $c_{i-1}$ is an in-twin of $c_{i+1}$ or not, we add $c_{i-1}$ to $\mathcal{B}$, resolving $v$ and $c_{i+k}$. Note that $c_{i-1}$ could have been added to $\mathcal{B}$ at the in-twin step, but we ensure that it is added in order to resolve the two conflicting vertices.
	
	In the fifth case, $\overrightarrow{C}$ contains exactly $2k$ vertices, and both $c_{i+k-1}$ and $c_{i+k+1}$ have out-neighbors $v_-$ and $v_+$, respectively, with in-degree~1. Assume furthermore that no vertex in $\overrightarrow{C} \setminus \{c_i\}$ has an in-arc, and that no vertex in $\overrightarrow{C} \setminus \{c_i, c_{i-1}, c_{i+1}\}$ is an in-twin. Now, if $c_{i-1}$ and $c_{i+1}$ are in-twins and have another in-twin, then, they will both be added to $\mathcal{B}$, and the pairs $(v_-,c_{i+k})$ and $(v_+,c_{i+k})$ are resolved. Otherwise, at least one of $v_-$ and $v_+$ will remain non-resolved with $c_{i+k}$, and thus we add $c_{i+k}$ to $\mathcal{B}$ in order to resolve the two pairs with the addition of just one vertex.
	
	Finally, for the sixth case, assume that $\overrightarrow{C}$ is of even length and contains two sources $c_i$ and $c_{i+2}$, that the two sinks are at equal distance of the sources (so they are $c_{i+1}$ and $c_{i+1+\frac{n}{2}}$), and that there is neither in-twin nor in-arc coming from outside of $\overrightarrow{C}$ (except the sources themselves). Now, if there are \annPaths, then, their endpoints and $c_{i+1+\frac{n}{2}}$ might not be resolved, since they are all at distance $\frac{n}{2}-1$ from the sources. Hence, we have to pick carefully among the potential sets of in-twins, and we may need to add another vertex to $\mathcal{B}$. There are three cases to cover.
	First, if all the \annPaths are \notSoBadPaths, then, by prioritizing the in-twins that are in the \notSoBadPaths, the endpoints and $c_{i+1+\frac{n}{2}}$ will be resolved without having to add another vertex in $\mathcal{B}$. Now, if all the \annPaths are \badPaths, then, every in-twin along the \annPaths belongs to a \annPath, and thus we will need to add a supplementary vertex to $\mathcal{B}$ in order to resolve the endpoints and the sink $c_{i+1+\frac{n}{2}}$; here, we choose the sink $c_{i+1}$. Finally, if there are both \badPaths and \notSoBadPaths, then, we can resolve the endpoints and the second sink by having the following priority on in-twins: those on \badPaths followed by those on \notSoBadPaths followed by those on non-\annPaths. Doing this will guarantee that the endpoints and $c_{i+1+\frac{n}{2}}$ are resolved. This settles the sixth and last special case.
	
	We will now prove that the vertices we added to $\mathcal{B}$ (sources, in-twins, and the six special cases) do form a resolving set. Since they were necessary to add, this will prove that $\mathcal{B}$ is indeed a metric basis.
	
	First, note that every vertex in $\overrightarrow{U}$ is reachable from some vertex in $\mathcal{B}$. Furthermore, recall that \Cref{alg-MDUnicyclic} gives priority to vertices in the cycle when resolving a set of in-twins, which will be important in some parts of the proof: if we know that the cycle contains an in-twin, we know that it will be in $\mathcal{B}$.
	
	Now, assume by contradiction that two vertices $u$ and $v$ are not resolved by $\mathcal{B}$. Since they are reachable, there is a vertex $b \in \mathcal{B}$ such that there are paths $P^b_u$ and $P^b_v$ from $b$ to $u$ and $v$, respectively. Let $x$ be the last common vertex of $P^b_u$ and $P^b_v$ (we can have $x=b$), and we can assume that every pair of predecessors of $u$ and $v$ is resolved (since, otherwise, we can just take the first unresolved pair of vertices on both paths). There are two cases to consider:
	\begin{enumerate}
		\item $u,v \in N^+(x)$. Since $u$ and $v$ are not resolved, they cannot be in-twins (since, otherwise, one of them would be in $\mathcal{B}$, a contradiction), and hence there is a vertex $w$ such that, without loss of generality, the arc $\overrightarrow{wu}$ exists and the arc $\overrightarrow{wv}$ does not exist. Now, $w$ has to be reachable from a vertex in $\mathcal{B}$, so there are two more possibilities.
		\\First, assume that there is a path from $x$ to $w$. There are two subcases here. In the first subcase, $u$, $w$ and $x$ are in the cycle $\overrightarrow{C}$ (of which $x$ is the only source and $u$ is the only sink). Then, whether $v$ is also in the cycle or not, either there is an in-twin or an in-arc coming from outside of $\overrightarrow{C}$ along the cycle, which would resolve $u$ and $v$, a contradiction; or we are in the third special case considered in \Cref{alg-MDUnicyclic}, and thus $u$ and $v$ are resolved, a contradiction. In the second subcase, the path from $x$ to $w$ goes through $u$, and thus $u$ and $w$ are both in the cycle $\overrightarrow{C}$ (which has no sink). Then, either there is an in-arc reaching $v$ or a vertex in $\overrightarrow{C}$, or there is an in-twin in $\overrightarrow{C}$, or we are in the second special case considered in \Cref{alg-MDUnicyclic}, and thus $u$ and $v$ are resolved, a contradiction.
		\\Now, assume that there is no such path, and thus there exists a vertex $b' \in \mathcal{B}$ such that there is a path of length $k$ from $b'$ to $w$. Since $u$ is an out-neighbor of $w$, this implies that there is a path of length $k+1$ from $b'$ to $v$. There are two subcases here. First, $x$ and $b'$ (or a representative along the path from $b'$ to $w$) are the two sources of the cycle $\overrightarrow{C}$, $u$ and $v$ are its two sinks, and $\overrightarrow{C}$ contains at least~6 vertices. In this subcase, it is necessary that, in $\overrightarrow{C}$, a vertex outside of the two sources is an in-twin or has an in-arc coming from outside of $\overrightarrow{C}$, which resolves $u$ and $v$, a contradiction. The second subcase is if the  path from $b'$ to $v$ goes through $x$ (whether $b=b'$ or not). But then, either there is an in-twin or an in-arc which resolves $u$ and $v$, a contradiction, or we are in the fourth or fifth special case considered in \Cref{alg-MDUnicyclic}, and thus $u$ and $v$ are resolved, a contradiction.
		
		\item $u,v \not\in N^+(x)$, hence, $u$ and $v$ each have a predecessor (respectively, $u'$ and $v'$) on $P^b_u$ and $P^b_v$. By hypothesis, there exists a vertex $b' \in \mathcal{B}$ that resolves $u'$ and $v'$ but not $u$ and $v$, so there is a path of length $k$ from $b'$ to $u'$ and a path of length $k+1$ from $b'$ to $v$ that does not go through $v'$. Note that $b'$ cannot be behind $x$ in $P^b_u$ or $P^b_v$ since otherwise it would not resolve $u'$ and $v'$, and it cannot be after in $P^b_u$ or $P^b_v$ since otherwise it would resolve $u$ and $v$, a contradiction. Hence, $b'$ (or a representative along the path from $b'$ to $v$) and $x$ are the two sources of the cycle $\overrightarrow{C}$, $u'$ and $v$ are its two sinks, and $\overrightarrow{C}$ contains at least six vertices. Like in the previous case, it is necessary that, either we are in the sixth special case considered in \Cref{alg-MDUnicyclic}, or, in $\overrightarrow{C}$, a vertex outside of the two sources is an in-twin or has an in-arc coming from outside of $\overrightarrow{C}$, which resolves $u$ and $v$, a contradiction.
	\end{enumerate}
	
	Hence, \Cref{alg-MDUnicyclic} resolves every pair of vertices in $\overrightarrow{U}$, and thus it returns a metric basis of a unicyclic graph. Finally, it is easy to see that it computes $\mathcal{B}$ in linear time; for the \annPaths in the sixth special case, we can do a breadth-first search of the graph starting from the sink $c_{i+1}$ to identify them, then go through the search tree again to compute, for each set of in-twins, which are all comprised of vertices in \annPaths (giving us \badPaths) and which are not (giving us \notSoBadPaths), and a third loop to correctly relabel the \annPaths.
\end{proof}

\section{FPT algorithm for modular width}

In a digraph $G$, a set $X \subseteq V(G)$ is a \emph{module} if every vertex not in $X$ 'sees' all vertices of $X$ in the same way. More precisely, for each $v \in V(G) \setminus X$ one of the following holds: (i) $(v,x),(x,v) \in E(G)$ for all $x \in X$, (ii) $(v,x),(x,v) \notin E(G)$ for all $x \in X$, (iii) $(v,x) \in E(G)$ and $(x,v) \notin E(G)$ for all $x \in X$, (iv) $(v,x) \notin E(G)$ and $(x,v) \in E(G)$ for all $x \in X$.
The singleton sets, $\emptyset$, and $V(G)$ are trivially modules of $G$. 
We call the singleton sets the \emph{trivial modules} of $G$. 

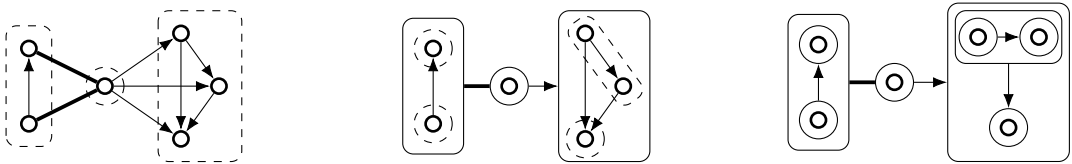
\begin{figure}[h]
	\centering
	\begin{subfigure}[b]{0.325\linewidth}
		\centering
		\begin{tikzpicture}
			\node[noeud] (c) at (0,0) {};
			\node[noeud] (a) at (-1,-.5) {};
			\node[noeud] (b) at (-1,.5) {};
			\node[noeud] (1) at (1,.7) {};
			\node[noeud] (2) at (1.5,0) {};
			\node[noeud] (3) at (1,-.7) {};
			\draw[\nicearrow] (a) to (b);
			\draw[line width=0.5mm] (a) to (c) to (b);
			\draw[\nicearrow] (c) to (1);
			\draw[\nicearrow] (c) to (2);
			\draw[\nicearrow] (c) to (3);
			\draw[\nicearrow] (1) to (2);
			\draw[\nicearrow] (1) to (3);
			\draw[\nicearrow] (2) to (3);
			\draw[dashed] (c) circle (0.25); 
			\draw[dashed,rounded corners] (-1.3,-.8) rectangle (-.7,.8); 
			\draw[dashed,rounded corners] (.7,1) rectangle (1.8,-1); 
		\end{tikzpicture}
	\end{subfigure}
	\begin{subfigure}[b]{0.325\linewidth}
		\centering
		\begin{tikzpicture}
			\node[noeud] (c) at (0,0) {};
			\node[noeud] (a) at (-1,-.5) {};
			\node[noeud] (b) at (-1,.5) {};
			\node[noeud] (1) at (1,.7) {};
			\node[noeud] (2) at (1.5,0) {};
			\node[noeud] (3) at (1,-.7) {};
			\node[draw, rectangle, rounded corners, minimum width=.8cm, minimum height=1.8cm] (X1) at (-1,0) {};
			\node[draw, circle, minimum width=.5cm] (X2) at (0,0) {};
			\node[draw, rectangle, rounded corners, minimum width=1.2cm, minimum height=2cm] (X3) at (1.25,0) {};
			\draw[\nicearrow] (a) to (b);
			\draw[line width=0.5mm] (X2) to (X1);
			\draw[\nicearrow] (X2) to (X3);
			\draw[\nicearrow] (1) to (2);
			\draw[\nicearrow] (1) to (3);
			\draw[\nicearrow] (2) to (3);
			\draw[dashed] (a) circle (0.25); 
			\draw[dashed] (b) circle (0.25); 
			\draw[dashed] (3) circle (0.25); 
			\draw[dashed, rounded corners, rotate around={-54:(1)}] (.8,.9) rectangle (2.1,.5); 
		\end{tikzpicture}
	\end{subfigure}
	\begin{subfigure}[b]{0.325\linewidth}
		\centering
		\begin{tikzpicture}
			\node[noeud] (c) at (0,0) {};
			\node[noeud] (a) at (-1,-.5) {};
			\node[noeud] (b) at (-1,.5) {};
			\node[noeud] (1) at (1.1,.6) {};
			\node[noeud] (2) at (1.9,.6) {};
			\node[noeud] (3) at (1.5,-.6) {};
			\node[draw, rectangle, rounded corners, minimum width=.8cm, minimum height=1.8cm] (X1) at (-1,0) {};
			\node[draw, circle, minimum width=.5cm] (X2) at (0,0) {};
			\node[draw, rectangle, rounded corners, minimum width=1.6cm, minimum height=2.1cm] (X3) at (1.5,0) {};
			\node[draw, circle, minimum width=.5cm] (Ma) at (a) {};
			\node[draw, circle, minimum width=.5cm] (Mb) at (b) {};
			\node[draw, circle, minimum width=.5cm] (M3) at (3) {};
			\node[draw, circle, minimum width=.5cm] (M1) at (1) {};
			\node[draw, circle, minimum width=.5cm] (M2) at (2) {};
			\node[draw, rectangle, rounded corners, minimum width=1.4cm, minimum height=.7cm] (X4) at (1.5,0.6) {};
			\draw[\nicearrow] (Ma) to (Mb);
			\draw[\nicearrow] (M1) to (M2);
			\draw[line width=0.5mm] (X2) to (X1);
			\draw[\nicearrow] (X2) to (X3);
			\draw[\nicearrow] (X4) to (M3);
		\end{tikzpicture}
	\end{subfigure}
	\caption{An example on how to decompose and draw a digraph using modules.}
	\label{fig-modularDecomposition}
\end{figure}

The graph $G[X]$ where $X$ is a module of $G$ is called a \emph{factor} of $G$. A family $\mathcal{X} = \{X_1, \ldots, X_s\}$ is a \emph{factorization} of $G$ if $\mathcal{X}$ is a partition of $V(G)$, and each $X_i$ is a module of $G$. If $X$ and $Y$ are two non-intersecting modules, then the relationship between $x \in X$ and $y \in Y$ is one of (i)-(iv) and always the same no matter which vertices $x$ and $y$ are exactly. Thus, given a factorization $\mathcal{X}$, we can identify each module with a vertex, and connect them to each other according to the arcs between the modules. More formally, we define the \emph{quotient} $G / \mathcal{X}$ with respect to the factorization $\mathcal{X}$ as the graph with the vertex set $\mathcal{X} = \{X_1, \ldots, X_s\}$ and $(X_i,X_j) \in E(G / \mathcal{X})$ if and only if $(x_i,x_j) \in E(G)$ where $x_i \in X_i$ and $x_j \in X_j$. A quotient depicts the connections of the different modules of a factorization to each other while omitting the internal structure of the factors. Each factor itself can be factorized further (as long as it is nontrivial, i.e. not a single vertex). By factorizing the graph $G$ and its factors until no further factorization can be done, we obtain a \emph{modular decomposition} of $G$. An example of a modular decomposition of a digraph is depicted on \Cref{fig-modularDecomposition}. The \emph{width} of a decomposition is the maximum number of sets in a factorization (or equivalently, the maximum number of vertices in a quotient) in the decomposition. The \emph{modular width} of $G$ is defined as the minimum width over all possible modular decompositions of $G$, and we denote it by $\mw (G)$. An optimal modular decomposition of a digraph can be computed in linear time~\cite{McConnellDirectedModular}. \MD for undirected graphs was shown to be fixed parameter tractable when parameterized by modular width by Belmonte et al.~\cite{MDwidth}. We will generalize their algorithm to directed graphs and strong and weak metric dimensions.

The following result lists several useful observations. 
\begin{proposition}\label{obs-modules}
	Let $\mathcal{X} = \{X_1, \ldots , X_s\}$ be a factorization of $G$, and let $W \subseteq V(G)$ be a resolving set of $G$. Denote $W_i = W \cap X_i$ for each $i \in \{1,\ldots,s\}$.
	\begin{enumerate}[label=(\roman*)]
		\item For all $x,y \in X_i$ and $z \in X_j$, $i\neq j$, we have $\dist_G(x,z) = \dist_G(y,z)$ and $\dist_G(z,x) = \dist_G(z,y)$.
		\item For all $x \in X_i$ and $y \in X_j$, $i\neq j$, we have $\dist_G(x,y) = \dist_{G/\mathcal{X}}(X_i,X_j)$.
		\item For all $x,y \in V(G)$ we have either $\dist_G(x,y) \leq \mw (G)$ or $\dist_G(x,y) = \infty$.
		\item The set $\{X_i \in \mathcal{X} \, | \, W_i \neq \emptyset \}$ is a resolving set of the quotient $G / \mathcal{X}$.
		\item For all distinct $x,y \in X_i$, where $X_i \in \mathcal{X}$ is non-trivial, we have $\dist_G(w,x) \neq \dist_G(w,y)$ for some $w \in W_i$.
		\item Let $w_1,w_2 \in X_i$. If $\dist_G(w_1,x) \neq \dist_G(w_2,x)$, then $x \in X_i$ and $\dist_G(w_1,x) \neq \dist_G(w_1,y)$ or $\dist_G(w_2,x) \neq \dist_G(w_2,y)$ for each $y \notin X_i$. 
	\end{enumerate}
\end{proposition}

The basic idea of our algorithm (and that of~\cite{MDwidth}) is to compute metric bases that satisfy certain conditions for the factors and combine these local solutions into a global solution. 
We know that non-trivial modules must contain elements of a resolving set, as modules must be resolved locally (\Cref{obs-modules}~(i)).
While combining the local solutions of non-trivial modules, we need to make sure that a vertex $x \in X_i$, where $X_i$ is non-trivial, is resolved from all $y \notin X_i$.
If $x$ and $y$ are resolved as described in~\Cref{obs-modules}~(vi), then we need to do nothing special.
However, if $x \in X_i$ is such that $\dist_G(w,x) = d$ for all $w \in W_i$ and a fixed $d \in \{1, \ldots, \mw (G), \infty\}$, there might exist a vertex $y \notin X_i$ such that $W_i$ does not resolve $x$ and $y$. 
We call such a vertex $x$ \emph{$d$-constant} (with respect to $W_i$).
We need to keep track of $d$-constant vertices and make sure they are resolved when we combine the local solutions.
There are at most $\mw(G) + 1$ $d$-constant vertices in each factor due to~\Cref{obs-modules}~(iii). 
We need to also make sure vertices in different modules that contain no elements of the solution set are resolved.
To do this, we might need to include some vertices from the trivial modules in addition to the vertices we have included from the non-trivial modules.

In the algorithm presented in~\cite{MDwidth}, the problems described above are dealt with by computing values $w(H,p,q)$ for every factor $H$, where $w(H,p,q)$ is the minimum cardinality of a resolving set of $H$ (with respect to the distance in $G$) where some vertex is 1-constant iff $p = true$ and some vertex is 2-constant iff $q = true$ (for undirected graphs these are the only two relevant cases). 
The same values are then computed for the larger graph by combining different solutions of the factors and taking their minimum. 
Our generalization of this algorithm is along the same lines as the original, however, we have more boolean values to keep track of. 
One difference to the techniques of the original algorithm is that we do not use the auxiliary graphs Belmonte et al. use. 
These auxiliary graphs were needed to simulate the distances of the vertices of a factor in $G$ as opposed to only within the factor. 
In our approach, we simply use the distances in $G$ and not the distances in the factors or the auxiliary graphs.

\begin{theorem}
	The metric dimension of a digraph $G$ with $\mw(G) \leq t$ can be computed in time $\mathcal{O}(t^5 2^{t^2} n + n^3 + m )$ where $n=|V(G)|$ and $m = |E(G)|$.
\end{theorem}
\begin{proof}
	Let us consider one level of an optimal modular decomposition of $G$. Let $H$ be a factor somewhere in the decomposition, and let $\mathcal{X} = \{X_1, \ldots , X_s\}$ be the factorization of $H$ according to the modular decomposition. For the graph $H$ (and its non-trivial factors $H[X_i]$) we denote by $w(H,\vp)$ the minimum cardinality of a set $W \subseteq V(H)$ such that
	\begin{enumerate}[label=(\roman*)]
		\item $W$ resolves $V(H)$ in $G$,
		\item $\vp = (p_1,\ldots,p_{\mw(G)},p_\infty)$ where $p_d = true$ if and only if $H$ contains a $d$-constant vertex with respect to $W$.
	\end{enumerate}
	If such a set does not exist, then $w(H,\vp)=\infty$. In order to compute the values $w(H,\vp)$, we next introduce the auxiliary values $\omega (\vp,I,P)$. The values $w(H[X_i],\vp)$ are assumed to be known for all $\vp$ and non-trivial modules $X_i$. Let the factorization $\mathcal{X}$ be labeled so that the modules $X_i$ are trivial for $i \in \{1, \ldots, h\}$ and non-trivial for $i \in \{h+1, \ldots , s\}$. Let $I \subseteq \{1, \ldots, h\}$ and 
	\[
	P =
	\left( {\begin{array}{c}
			\vp^{h+1} \\
			\vdots  \\
			\vp^{s}  \\
	\end{array} } \right).
	\]
	We define $\omega (\vp,I,P) = |I| + \sum_{i=h+1}^{s} w(H[X_i],\vp^i)$ if the following conditions (a)-(d) hold. In what follows, a representative of a module $X_i$ is denoted by $x_i$.
	\begin{enumerate}[label=(\alph*)]
		\item The set $Z=\{X_i \in \mathcal{X} \, | \, i \in I \cup \{h+1,\ldots,s\} \}$ resolves the quotient $H/\mathcal{X}$ with respect to the distances in $G$.
		\item For $d \in \{1,\ldots,\mw(G),\infty\}$ and $i \in \{h+1,\ldots,s\}$, if $p_d^i = true$, then for each trivial module $X_j = \{x_j\}$ where $j \notin I$ we have $\dist_G(x_i,x_j) \neq d$ or there exists $X_k \in Z \setminus \{X_i\}$ such that $\dist_G(x_k,x_i) \neq \dist_G(x_k,x_j)$.
		\item For $d_1,d_2 \in \{1,\ldots,\mw(G),\infty\}$ and distinct $i,j \in \{h+1,\ldots,s\}$, if $p_{d_1}^i = p_{d_2}^j = true$, then $\dist_G(x_i,x_j) \neq d_1$, or $\dist_G(x_j,x_i) \neq d_2$, or there exists $X_k \in Z \setminus \{X_i, X_j\}$ such that $\dist_G(x_k,x_i) \neq \dist_G(x_k,x_j)$.
		\item For all $d \in \{1,\ldots,\mw(G),\infty\}$, we have $p_d = true$ (in $\vp$) if and only if for some $i \in \{1, \ldots , h\} \setminus I$ we have $\dist_G(x_j,x_i) = d$ for all $X_j \in Z$, or for some $i \in \{h+1, \ldots , s\}$ we have $p_d^i= true$ and $\dist_G(x_j,x_i) = d$ for all $X_j \in Z \setminus \{X_i\}$.
	\end{enumerate}
	If these conditions cannot be met, then we set $\omega (\vp,I,P) = \infty$.
	
	Here we give an outline of a proof for the equality $w(H,\vp) = \min_{I,P} \omega (\vp,I,P)$. 
	
	We will first show that $w(H,\vp) \leq \min_{I,P} \omega (\vp,I,P)$. This clearly holds if $\min_{I,P} \omega (\vp,I,P) = \infty$. 
	So assume that $I$ and $P$ are such that $\omega (\vp,I,P)$ is as small as possible. Let $W \subseteq V(G)$ be such that each $W_i = X_i \cap W$ is a $w(H[X_i],\vp^i)$-set for $i \in \{h+1,\ldots,s\}$, $|W| = \omega (\vp,I,P)$, and $W$ fulfills the conditions (a)-(d). We will show that $W$ fulfills the conditions (i) and (ii), and thus $|W| \geq w(H,\vp)$.
	\begin{enumerate}[label=(\roman*)]
		\item $W$ resolves $V(H)$ in $G$: If $x,y \in X_i$, then $i \in \{h+1,\ldots,s\}$ and $x$ and $y$ are resolved by $W_i$. Assume that $x \in X_i$ and $y \in X_j$, $i \neq j$. Suppose that $x$ and $y$ are not resolved due to~\Cref{obs-modules}~(vi). Then at least one of them is $d$-constant or they are both in trivial modules. Now, if $i,j \in \{1,\ldots,h\}$, then $x$ and $y$ are resolved due to condition (a). If $i \in \{1, \ldots ,h\}$ and $j \in \{h+1,\ldots, s\}$, then $x$ and $y$ are resolved due to condition (b). If $i,j \in \{h+1,\ldots, s\}$, then $x$ and $y$ are resolved due to condition (c).
		\item This holds due to condition (d).
	\end{enumerate}
	Therefore, $w(H,\vp) \leq \min_{I,P} \omega (\vp,I,P)$.
	
	Let us then show that $w(H,\vp) \geq \min_{I,P} \omega (\vp,I,P)$. Again, if $w(H,\vp) = \infty$, then the claim clearly holds. So assume that $W \subseteq V(G)$ is such that $|W|=w(H,\vp)$ and $W$ fulfills conditions (i) and (ii). Consider the sets $W_i$ for $i \in \{h+1,\ldots,s\}$. Each $W_i$ resolves $X_i$ in $G$ due to~\Cref{obs-modules}~(v), and thus $|W_i| \geq w(H[X_i],\vp^i)$ for $\vp^i$ defined with respect to $W_i$. Let $P$ be defined with these $\vp^i$'s, and let $I = \{i \in \{1, \ldots , h\} \, | \, W_i \neq \emptyset\}$. Now, $|W| \geq |I| + \sum_{i=h+1}^{s} w(H[X_i],\vp^i)$. Moreover, the conditions (a)-(d) hold:
	\begin{enumerate}[label=(\alph*)]
		\item Holds due to condition (i).
		\item Assume to the contrary that $p_d^i = true$ and $j\in \{1, \ldots, h\} \setminus I$ is such that $\dist_G(x_k,x_i) = \dist_G (x_k,x_j)$ for all $X_k \in Z \setminus \{X_i\}$. Let $x \in X_i$ be $d$-constant. Since $W$ resolves $x$ and $y$, there exists $w \in W_i$ such that $\dist_G(w,x) \neq \dist_G(w,x_j) $, and thus $\dist_G(x_i,x_j) = \dist_G(w,x_j) \neq d$.
		\item Can be shown with the same technique as (b).
		\item Clear.
	\end{enumerate}
	Therefore, $|W| \geq \omega (\vp,I,P)$ and $w(H,\vp) \geq \min_{I,P} \omega (\vp,I,P)$.
	
	Let us then discuss the complexity of this algorithm. As a preprocessing step, we need to compute the distances between all pairs of vertices. This can be done using the Floyd-Warshall algorithm in $\mathcal{O}(n^3)$ time. An optimal modular decomposition can be computed in $\mathcal{O}(n+m)$ time~\cite{McConnellDirectedModular}. We then need to compute the values $w(H,\vp)$ for each factor $H$ starting from the trivial modules and working our way up in the decomposition. The values $w(H,\vp)$ are computed using the auxiliary values $\omega (\vp,I,P)$. There are $\mathcal{O}(2^{t^2})$ different possibilities for $I$ and $P$ and their combinations (note that the vector $\vp$ is determined based on $I$ and $P$). For each pair $I$,$P$, we need to check the conditions (a)-(d), out of which (c) is the most costly time-wise and can be checked in $\mathcal{O}(t^5)$ time. Thus, computing the values $w(H,\vp)$ can be done in $\mathcal{O}(t^5 2^{t^2})$ time for each $H$. The total computing time then follows from the fact that there are at most $2n$ factors in any modular decomposition of a graph with $n$ vertices. (The decomposition can be presented as a rooted tree where the vertices represent the factors and edges represent inclusion. In this tree the leaves are exactly the trivial modules, and there are $n$ of them. Every internal vertex has degree at least 3, except the root has degree at least 2. Using the handshake lemma it is then straightforward to show that this tree can have at most $2n$ vertices.)
\end{proof}

The original algorithm of Belmonte et al. has conditions (a)-(g), of which (a) is (essentially) the same as (a) above, (b) and (c) are covered by (b), (d) and (e) by (c), and (f) and (g) by (d). Notice that our condition (c) is true whenever $\dist_G(x_i,x_j) = \dist_G(x_j,x_i)$ and $d_1 \neq d_2$. Specifically, if $G$ is undirected, we do not need to care about (c) for pairs where $d_1 \neq d_2$.

\section{NP-hardness for restricted DAGs}

We now complement the hardness result from~\cite{ABCHONS20}, which was for bipartite DAGs of maximum degree~8 and maximum distance~4.

\begin{theorem}
	\MD is NP-complete, even on planar triangle-free DAGs of maximum degree~6 and maximum distance~4.
\end{theorem}

\begin{proof}
	The problem is clearly in NP: a certificate is a set of vertices, for which we can check in polynomial time if it is of the required size and if it resolves all vertices by computing the distance vectors and comparing them.
	
	For completeness, we reduce from \textsc{Vertex Cover} on 2-connected planar cubic graphs, which is known to be NP-complete~\cite[Theorem 4.1]{M01}.
	
	Given a 2-connected planar cubic graph $G$, we construct a DAG $G'$ as follows. First of all, note that by Petersen's theorem, $G$ contains a perfect matching $M\subset E(G)$, that can be constructed in polynomial time. A planar embedding of $G$ can also be constructed in polynomial time, so we fix one. We let $V(G')=V(G)\bigcup_{e=uv\in E(G)}\{a_e,b_e,c_e,d_e^u,d_e^v\}\bigcup_{e=uv\in M}\{f_e,g_e,h_e\}$. For every edge $e=uv$ of $G$, we add the arcs $\{\overrightarrow{a_eb_e},\overrightarrow{b_ec_e},\overrightarrow{c_ed_e^u},\overrightarrow{c_ed_e^v},\overrightarrow{ud_e^u},\overrightarrow{vd_e^v}\}$. For every edge $e=uv$ of the perfect matching $M$ of $G$, assuming the neighbors of $u$ (in the clockwise cyclic order with respect to the planar embedding of $G$) are $v,x,y$ and those of $v$ are $u,s,t$, we arbitrarily fix one side of the edge $uv$ to place the vertices $f_e$, $g_e$ and $h_e$ (say, on the side that is close to the edges $ux$ and $vt$). We add the arcs $\{\overrightarrow{f_eg_e},\overrightarrow{g_ec_e},\overrightarrow{g_eh_e},\overrightarrow{h_eu},\overrightarrow{h_ev},\overrightarrow{c_ec_{uy}},\overrightarrow{c_ec_{vs}},\overrightarrow{h_ec_{ux}},\overrightarrow{h_ec_{vt}}\}$.
	
	Using the embedding of $G$, $G'$ can also be drawn in a planar way, it has maximum degree~6 (the vertices of type $c_e$ are of degree~6 when $e\in M$), has no triangles, and no shortest directed path of length~5. See Figure~\ref{fig:reduc} for an illustration.
	
	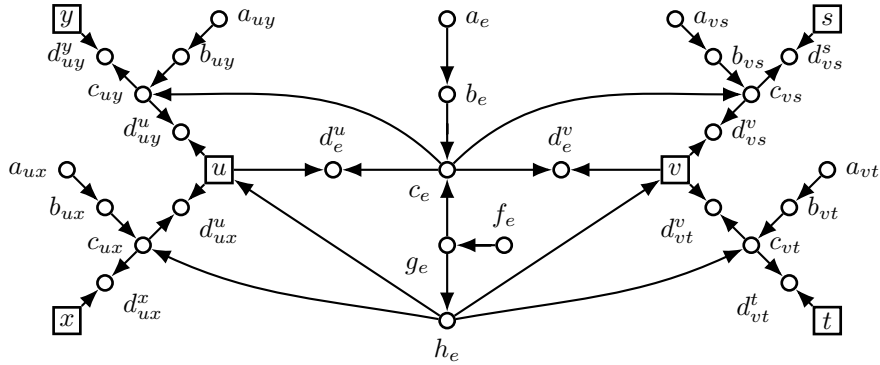
\begin{figure}[!htpb]
		\centering
		\begin{tikzpicture}
			
			\node[bignoeud,shape=rectangle](u) at (-3,0)    {$u$};
			\node[bignoeud,shape=rectangle](v) at (3,0)    {$v$};
			
			\node[noeud,label={90:$d_e^u$}](du) at (-1.5,0)    {};
			\node[noeud,label={90:$d_e^v$}](dv) at (1.5,0)    {};	
			\node[noeud,label={0:$a_e$}](a) at (0,2)    {};		
			\node[noeud,label={0:$b_e$}](b) at (0,1)    {};		
			\node[noeud,label={225:$c_e$}](c) at (0,0)    {};		
			\node[noeud,label={90:$f_e$}](f) at (0.75,-1)    {};		
			\node[noeud,label={200:$g_e$}](g) at (0,-1)    {};		
			\node[noeud,label={270:$h_e$}](h) at (0,-2)    {};
			
			\path[thick,\nicearrow] (u) edge[out=0,in=180,min distance=10mm] node {} (du);
			\path[thick,\nicearrow] (c) edge[out=0,in=180,min distance=10mm] node {} (dv);
			\path[thick,\nicearrow] (v) edge[out=180,in=0,min distance=10mm] node {} (dv);
			\path[thick,\nicearrow] (c) edge[out=180,in=0,min distance=10mm] node {} (du);
			\path[thick,\nicearrow] (a) edge node {} (b);
			\path[thick,\nicearrow] (b) edge[out=270,in=90,min distance=10mm] node {} (c);
			\path[thick,\nicearrow] (f) edge[out=180,in=0,min distance=10mm] node {} (g);
			\path[thick,\nicearrow] (g) edge node {} (c);
			\path[thick,\nicearrow] (g) edge node {} (h);
			\path[thick,\nicearrow] (h) edge node {} (u);
			\path[thick,\nicearrow] (h) edge node {} (v);
			
			\node[bignoeud,shape=rectangle](y) at (-5,2)    {$y$};
			\node[noeud,label={180:$d_{uy}^u$}](duyu) at (-3.5,0.5)    {};
			\node[noeud,label={180:$c_{uy}$}](cuy) at (-4,1)    {};
			\node[noeud,label={180:$d_{uy}^y$}](duyy) at (-4.5,1.5)    {};
			\node[noeud,label={0:$a_{uy}$}](auy) at (-3,2)    {};		
			\node[noeud,label={0:$b_{uy}$}](buy) at (-3.5,1.5)    {};	
			
			\path[thick,\nicearrow] (u) edge node {} (duyu);
			\path[thick,\nicearrow] (cuy) edge node {} (duyu);
			\path[thick,\nicearrow] (cuy) edge node {} (duyy);
			\path[thick,\nicearrow] (auy) edge node {} (buy);
			\path[thick,\nicearrow] (buy) edge node {} (cuy);
			\path[thick,\nicearrow] (y) edge node {} (duyy);

			\node[bignoeud,shape=rectangle](x) at (-5,-2)    {$x$};
			\node[noeud,label={-10:$d_{ux}^u$}](duxu) at (-3.5,-0.5)    {};
			\node[noeud,label={180:$c_{ux}$}](cux) at (-4,-1)    {};
			\node[noeud,label={-10:$d_{ux}^x$}](duxx) at (-4.5,-1.5)    {};	
			\node[noeud,label={180:$a_{ux}$}](aux) at (-5,0)    {};		
			\node[noeud,label={180:$b_{ux}$}](bux) at (-4.5,-0.5)    {};	
			
			\path[thick,\nicearrow] (u) edge node {} (duxu);
			\path[thick,\nicearrow] (cux) edge node {} (duxu);
			\path[thick,\nicearrow] (cux) edge node {} (duxx);
			\path[thick,\nicearrow] (aux) edge node {} (bux);
			\path[thick,\nicearrow] (bux) edge node {} (cux);
			\path[thick,\nicearrow] (x) edge node {} (duxx);

			\node[bignoeud,shape=rectangle](s) at (5,2)    {$s$};
			\node[noeud,label={0:$d_{vs}^v$}](dvsv) at (3.5,0.5)    {};
			\node[noeud,label={0:$c_{vs}$}](cvs) at (4,1)    {};
			\node[noeud,label={0:$d_{vs}^s$}](dvss) at (4.5,1.5)    {};
			\node[noeud,label={0:$a_{vs}$}](avs) at (3,2)    {};		
			\node[noeud,label={0:$b_{vs}$}](bvs) at (3.5,1.5)    {};	
			
			\path[thick,\nicearrow] (v) edge node {} (dvsv);
			\path[thick,\nicearrow] (cvs) edge node {} (dvsv);
			\path[thick,\nicearrow] (cvs) edge node {} (dvss);
			\path[thick,\nicearrow] (s) edge node {} (dvss);	
			\path[thick,\nicearrow] (avs) edge node {} (bvs);
			\path[thick,\nicearrow] (bvs) edge node {} (cvs);	
			
			\node[bignoeud,shape=rectangle](t) at (5,-2)    {$t$};
			\node[noeud,label={190:$d_{vt}^v$}](dvtv) at (3.5,-0.5)    {};
			\node[noeud,label={0:$c_{vt}$}](cvt) at (4,-1)    {};
			\node[noeud,label={190:$d_{vt}^t$}](dvtt) at (4.5,-1.5)    {};
			\node[noeud,label={0:$a_{vt}$}](avt) at (5,0)    {};		
			\node[noeud,label={0:$b_{vt}$}](bvt) at (4.5,-0.5)    {};		
			
			\path[thick,\nicearrow] (v) edge node {} (dvtv);
			\path[thick,\nicearrow] (cvt) edge node {} (dvtv);
			\path[thick,\nicearrow] (cvt) edge node {} (dvtt);
			\path[thick,\nicearrow] (avt) edge node {} (bvt);
			\path[thick,\nicearrow] (bvt) edge node {} (cvt);
			\path[thick,\nicearrow] (t) edge node {} (dvtt);	
			
			\path[thick,\nicearrow] (c) edge[out=45,in=180] node {} (cvs); 
			\path[thick,\nicearrow] (c) edge[out=135,in=0] node {} (cuy); 		
			\path[thick,\nicearrow] (h) edge[out=10,in=205] node {} (cvt); 
			\path[thick,\nicearrow] (h) edge[out=170,in=-25] node {} (cux);  		
			
		\end{tikzpicture}\centering
		\caption{Illustration of the reduction for an edge $e=uv$ of $G$ in $M$, and the surrounding edges $ux$, $uy$, $vs$ and $vt$. Squared vertices are the original ones from $G$.}\label{fig:reduc}
	\end{figure}
	
	Now, we claim that $G$ has a vertex cover of size at most $k$ if and only if $G'$ has metric dimension at most $k+|E(G)|+|M|=k+4|E(G)|/3=k+2|V(G)|$.
	
	If $G$ has a vertex cover $C$ of size $k$, we construct a resolving set $R(C)$ of $G'$ as follows. Include the vertices of $C$ in $R(C)$, as well as all vertices of $\{a_e~|~e\in E(G)\}\cup\{f_e~|~e\in M\}$. The vertices in $R(C)$ are clearly uniquely resolved. For a given edge $e$ of $G$, the vertices $b_e$ and $c_e$ are uniquely at distance~1 and~2 from $a_e$, respectively, so all vertices of these types are uniquely resolved. Among the other vertices associated to $e$, $d_e^u$ and $d_e^v$ are the only ones at distance~3 from $a_e$; moreover, $d_e^u$ is at distance~1 from $u$ and $d_e^v$ at distance~1 from $v$, but not vice-versa, so $C\cap\{u,v\}$ resolves $d_e^u$ and $d_e^v$. Thus, all vertices of these types are uniquely resolved. Among the remaining vertices, $g_e$ and $h_e$ are uniquely at distance~1 and~2 from $f_e$, respectively, so all vertices of these types are uniquely resolved. Finally, the vertices in $V(G)\setminus R(C)$ are resolved by the unique vertex of type $f_e$ from which each of them is at distance~3. Hence, all vertices are uniquely resolved and $R(C)$ is indeed a resolving set of $G'$.
	
	Conversely, let $R$ be a resolving set of $G'$ of size at most $k+|E(G)|+|M|$. Notice that for each edge $e$ of $G$, one of $a_e,b_e$ belongs to $R$ in order to resolve this pair, and similarly, for each edge $e$ in $M$, one of $f_e,g_e$ belongs to $R$ (in the case of strong metric dimension, $a_e$ and $f_e$ belong to the solution, since they are sources).
	
	We construct a potential vertex cover $C(R)$ by taking $R\cap V(G)$. Moreover, for each edge $e=uv$ in $M$, we add to $C(R)$ any of $u,v$ (if possible, one that is not yet in $C(R)$) in case the set $\{a_e,b_e,c_e,f_e,g_e,h_e,d_e^u,d_e^v\}$ contains three vertices of $R$. If it contains at least four, both $u,v$ are put into $C(R)$. Similarly, for each edge $e=uv$ in $E(G)\setminus M$, we add to $C(R)$ any of $u,v$ in case the set $\{a_e,b_e,c_e,d_e^u,d_e^v\}$ contains two vertices of $R$ (if possible, we add one that is not yet in $C(R)$), and we add both $u,v$ if it contains more than two vertices of $R$.
	
	By the above paragraph, the resulting set $C(R)$ contains at most $|R|-|E(G)|-|M|\leq k$ vertices. Now, consider a pair $d_e^u,d_e^v$ for some edge $e$. If $u$ or $v$ is in $R$, it is also in $C(R)$, and $e$ is covered by $C(R)$. Assume now that none of $u,v$ is in $R$. If $e\in M$, necessarily one of $u,v,d_e^u,d_e^v$ belongs to $R$ to resolve that pair, and so, as none of $u,v$ are in $R$, $|a_e,b_e,c_e,f_e,g_e,h_e,d_e^u,d_e^v|\geq 3$ and by our construction, either $u$ or $v$ (or both) have been added to $C(R)$. Thus, $e$ is covered by $C(R)$. If $e\notin M$, the only vertices that can resolve $d_e^u,d_e^v$ are again $u,v,d_e^u,d_e^v$, or a vertex $h_{e'}$ where $e'\neq e$ is an edge of $G$ in $M$ incident with $u$ or $v$ and $h_{e'}$ is not adjacent to $c_{e}$. Again, as none of $u,v$ is in $R$, if one of $d_e^u,d_e^v$ belongs to $R$, $|a_e,b_e,c_e,d_e^u,d_e^v|\geq 2$ and by our construction, either $u$ or $v$ (or both) have been added to $C(R)$. Otherwise, it must be that some vertex $h_{e'}$ is in $R$, where $e'\neq e$ is an edge of $G$ in $M$ incident with $u$ or $v$ (say, $u$ and $e'=uw$) and $h_{e'}$ is not adjacent to $c_{e}$. But notice that $h_{e'}$ does not resolve $d_{e'}^u$ and $d_{e'}^w$, as $e'\in M$. Thus, either $w\in R$ and $|a_{e'},b_{e'},c_{'e},f_{e'},g_{e'},h_{e'},d_{e'}^u,d_{e'}^w|\geq 3$, or $w\notin R$ and $|a_{e'},b_{e'},c_{'e},f_{e'},g_{e'},h_{e'},d_{e'}^u,d_{e'}^w|\geq 4$. In both cases, by our construction, we would have added $u$ to $C(R)$. Thus, in all cases, one of $u,v$ belongs to $C(R)$ and $e$ is covered. Thus, $C(R)$ is a vertex cover of size at most $k$, as needed.
\end{proof}

\section{Conclusion}

\MD can be solved in polynomial time on outerplanar graphs, using an involved algorithm~\cite{MDplanar}. Can one generalize our algorithms for trees and unicyclic graphs to solve \MD for directed (or at least, oriented) outerplanar graphs in polynomial time? Extending our algorithm to cactus graphs already seems non-trivial.

One open question is whether \MD is NP-hard on planar bipartite subcubic DAGs?

Also, it would be interesting to see which hardness results known for \MD of undirected graphs also hold for DAGs, or for oriented graphs.

\section*{Acknowledgement}

We thank the anonymous referee for their careful reading of the proofs and their helpful suggestions.

\section*{Dedication to Rolf Niedermeier}

This work may not have been performed without the influence of Rolf Niedermeier on Florent Foucaud. As a PhD student in 2012, Florent visited Rolf's group in TU Berlin for two weeks, and again for a month in 2013. Despite Florent not being experienced in the field of parameterized complexity, Rolf warmly welcomed these visits and made Florent feel at ease. Florent started a collaboration on the parameterized complexity of the \MD problem with some members of Rolf's group (inspired by two of them, who had just obtained an important result in the area~\cite{HN13}). This collaboration did not lead to any publication, nevertheless, the discussions with Rolf, his students and visitors and the friendly atmosphere in the group certainly influenced Florent's later research and inspired him to work more on the parameterized complexity of graph problems. This includes the present paper, which is also about \MD. Rolf was also a positive model on how to have a dynamic and positive research group. We have been very saddened by Rolf's unexpected and too early passing, and dedicate this paper to his memory.

\bibliographystyle{abbrv}
\bibliography{references}

\end{document}